\PassOptionsToPackage{usenames,dvipsnames}{xcolor}
\documentclass[a4paper, parskip=half, twocolumn=off, twoside=on]{scrartcl}
\usepackage[hidelinks]{hyperref}
\usepackage{geometry}

\usepackage[english]{babel}
\usepackage{amsmath}
\usepackage{amsthm}
\usepackage[fixamsmath]{mathtools}
\mathtoolsset{showonlyrefs=true}
\usepackage{booktabs}%

\usepackage{orcidlink}
\usepackage{siunitx}
\usepackage{float}
\usepackage{caption}
\usepackage{xfrac}
\usepackage{booktabs}
\usepackage{xcolor}%
\usepackage[many]{tcolorbox}
\usepackage{csquotes}
\usepackage[backend=biber, style=numeric, sortcites=true, citestyle=numeric-comp, giveninits=true, sorting=nyt, doi=true, isbn=false, url=false, maxcitenames=3, maxbibnames=10, minalphanames=3, block=none, date=year]{biblatex}
\AtEveryBibitem{\clearlist{language}\clearfield{pages}\clearfield{pagetotal}\clearfield{file}\clearlist{address}\clearlist{location}\clearlist{venue}}
 
\renewbibmacro{in:}{}
\usepackage[proportional]{libertinus}
\usepackage{libertinust1math}
\usepackage[cal=boondox,scr=boondox, bb=ams]{mathalpha}

\newcommand{\symbf}[1]{\mathbfit{#1}}
\newcommand{\symsf}[1]{\mathsf{#1}}
\newcommand{\symbfsf}[1]{\mathsfbf{#1}}
\newcommand{\symbfup}[1]{\mathbf{#1}}
\newcommand{\symbb}[1]{\mathbb{#1}}

\newcommand{\symcal}[1]{\mathcal{#1}}
\usepackage{authblk}

\theoremstyle{definition}
\newtheorem{definition}{Definition}

\theoremstyle{plain}

\theoremstyle{remark}
\newtheorem{remark}[definition]{Remark}
\newtheorem{example}[definition]{Example}
\tcbuselibrary{theorems}
\tcbset{
breakable, 
enhanced jigsaw, 
lines before break=20,
colback=white,
colframe=black!15, 
coltitle=black, 
colbacktitle=black!15, 
fonttitle=\bfseries\sffamily, arc=0mm,
subtitle style={colback=black!15,fontupper={\sffamily\color{black}}}
}
\tcbset{shield externalize}

\newtcolorbox[use counter*=definition]{boxalgorithm}[2][]{
 list text={#2},
 title={Algorithm \thetcbcounter: #2},
 #1
}

\usepackage[clines, plainfootsepline, automark]{scrlayer-scrpage}
\clearpairofpagestyles

\ihead[]{}
\cehead{J. Magiera, C. Rohde}
\cohead{\csname @title\endcsname}
\ohead[]{\pagemark}
\ifoot[]{}
\ofoot[]{}

\setkomafont{pagehead}{\normalfont\footnotesize\itshape}
\setkomafont{pagefoot}{\normalfont\footnotesize}
\setkomafont{pagenumber}{\normalfont\normalsize\itshape}
\KOMAoptions{parskip=half}
\KOMAoption{parskip}{half}
\KOMAoptions{abstract=true}

\allowdisplaybreaks[1] 
\newcommand{\pliq}{-}
\newcommand{\pvap}{+}
\newcommand{\phasel}{\pliq}
\newcommand{\phaser}{\pvap}

\newcommand{\phaserl}{\pm}
\newcommand{\cres}{CRes} %

\newcommand{\frictionvariable}{\ensuremath{f}}

\newcommand{\mdtimevar}{\tau} 
\newcommand{\mdinterface}{\Gamma} 
 
\newcommand{\mdtimeend}{\mdtimevar_{\mathrm{end}}} 
\newcommand{\mdstepend}{N_{\mathrm{end}}} 

\newcommand{\continterface}{\symbfup{\Gamma}} 
\newcommand{\contdomain}{\symbfup{\Omega}} 

\newcommand{\nnparams}{\symbfup{\theta}} 
\newcommand{\datasetvar}{D}

\newcommand{\fractionvar}{\alpha} 

\newcommand{\networkF}{\symsf{R}}

\newcommand{\pcsaft}{{PC-SAFT}}
\newcommand{\cmpa}{0}
\newcommand{\cmpb}{1}

\newcommand{\setsep}{\allowbreak:\allowbreak}

\providecommand\mapsfrom{\mathrel{\reflectbox{\ensuremath{\mapsto}}}}

\newcommand{\rmdmix}{\ensuremath{R_{\mathrm{MD}}}} 
\newcommand{\genericsolversym}{\symcal{R}} 
\newcommand{\PP}{\mathcal{P}}

\newcommand{\uu}{\vec{u}} 
\newcommand{\UU}{\vec{U}} 

\renewcommand{\vec}[1]{\symbf{#1}}

\newcommand{\mat}[1]{\symbfsf{#1}}

\newcommand{\bN}{\symbb{N}}

\newcommand{\bR}{\symbb{R}}
\newcommand{\bS}{\symbb{S}}

\newcommand{\cP}{\mathcal{P}}

\newcommand{\cR}{\mathcal{R}}

\newcommand{\cT}{\mathcal{T}}

\newcommand{\va}{\vec{a}}

\newcommand{\vm}{\vec{m}}
\newcommand{\vn}{\vec{n}}

\providecommand{\vv}{}
\renewcommand{\vv}{\vec{v}}

\newcommand{\vx}{\vec{x}}

\DeclareTextSymbolDefault{\DH}{T1}
\newcommand{\mathDJ}{\ensuremath{\text{\itshape\DH}}}

\DeclarePairedDelimiter\abs{\lvert}{\rvert}

\DeclarePairedDelimiter\norm{\lVert}{\rVert}

\DeclarePairedDelimiterX{\mean}[1]{\{}{\}}{\mkern-3.0mu\delimsize\{ {#1} \delimsize\}\mkern-3.0mu}

 \DeclareSIUnit{\redenergy}{\ensuremath{\varepsilon}}
\DeclareSIUnit{\redlength}{\ensuremath{\sigma}}
 \DeclareSIUnit[]{\boltzmannunit}{\text{\ensuremath{k_{\textup{B}}}}}
 \DeclareSIUnit{\atomicmassunit}{u}

\makeatletter
\DeclareRobustCommand{\AA}{%
 \leavevmode
 \vbox{\ialign{##\cr
 \hidewidth\char'27 \hidewidth\cr
 \noalign{\nointerlineskip\kern-1.4ex}
 A\cr
 }}%
}
\makeatother
\DeclareSIUnit{\angstrom}{\textup{\AA}}
\definecolor{argoncolor}{RGB}{27,158,119}
\definecolor{methanecolor}{RGB}{117,112,179}
\begin{document}
\title{A Multiscale Method for Two-Component, Two-Phase Flow with a Neural Network Surrogate}
\author{Jim Magiera\,\orcidlink{0000-0001-9807-0784}}
\author{Christian Rohde\,\orcidlink{0000-0001-9183-5094}}
\affil{University of Stuttgart, Institute of Applied Analysis and Numerical Simulation}
\date{}

\maketitle

\begin{abstract}
    \noindent{}Understanding the dynamics of phase boundaries in fluids requires quantitative knowledge about the microscale processes at the interface.
    We consider the sharp-interface motion of compressible two-component flow, and propose a heterogeneous multiscale method (HMM) to describe the flow fields accurately.
    The multiscale approach combines a hyperbolic system of balance laws on the continuum scale with molecular-dynamics simulations on the microscale level.
    Notably, the multiscale approach is necessary to compute the interface dynamics because there is -- at present -- no closed continuum-scale model.
    The basic HMM relies on a moving-mesh finite-volume method, and has been introduced recently for compressible one-component flow with phase transitions in \cite{MagieraRohde22}.
    To overcome the numerical complexity of the molecular-dynamics microscale model a deep neural network is employed as an efficient surrogate model.\\
    The entire approach is finally applied to simulate droplet dynamics for argon--methane mixtures in several space-dimensions. Up to our knowledge such compressible two-phase dynamics accounting for microscale phase-change transfer rates have not yet been computed.
\end{abstract}
\section{Introduction}\label{intro}
\begin{figure}[htb]%
    \centering
    {\includegraphics[width=0.236076\columnwidth]{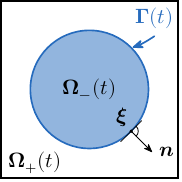}}
    \caption{Partitioned domain $\contdomain$ with interface $\boldsymbol{\Gamma}(t)$ and bulk domains
        $\contdomain_{\pm}$ for $d=2$ spatial dimensions.}
    \label{fig:geometry}
\end{figure}

In this work we consider the dynamics of a compressible two-component fluid which allows for phase changes between a liquid and a vapor phase state.
Keeping temperature constant, we aim to determine continuum-scale quantities like component-wise densities and momenta on a spatial scale where all phase boundaries are explicitly represented, i.e., without further
averaging like in e.g. Baer--Nunziato modelling \cite{saurel.abgrall:multiphase:1999,Andrianov2004b,Zein2010}.
This choice implies that the phase -- liquid or vapor -- is uniquely determined by the values of the component densities.
Having made this decision, we focus on a sharp-interface approach, and refer to \cite{AMW:98} and just recently \cite{KMR:22} for diffuse-interface models.
We start from an extension of the isothermal Euler equations with frictional forces for one-phase flows, as introduced in \cite{bothe.dreyer:continuum:2015}, to the two-phase case, see \eqref{eq:iso_euler_multicomponent} below.
We will focus mainly on the two-component system comprised of the noble gases argon and methane which allows a relatively simple thermodynamical description.
These evolution equations are supposed to hold in the liquid-phase domain \(\contdomain_{\pliq}(t)\) and the vapor-phase domain \(\contdomain_{\pvap}(t)\), where they form a set of hyperbolic balance laws.
The two bulk domains are separated by a co-dimension-1 manifold, the interface $\boldsymbol{\Gamma} = \boldsymbol{\Gamma}(t)$.
Thus, the entire fluid domain \(\contdomain \subset \bR^d\), \(d \in \{1,2,3\}\), at any time $t \in [0,t_{\mathrm{end}}]$, $t_{\mathrm{end}}>0 $, is partitioned into
$\contdomain = \contdomain_{\pvap}(t) \cup \boldsymbol{\Gamma}(t) \cup \contdomain_{\pvap}(t)$, see Figure~\ref{fig:geometry} for illustration.
The basic mathematical model is provided in Section~\ref{ch:multicomponent}.
\captionsetup{width=0.75\textwidth}
To determine the motion of the interface and to ensure the well-posedness of the resulting free boundary value problem one has to prescribe transmission conditions across the interface $\boldsymbol{\Gamma}(t)$.
Obviously, the Rankine--Hugoniot conditions have to be imposed; but it is well-known, even in the one-component case, that additional conditions have to be prescribed at the interface \cite{CP:03,Han2013,hantke.thein:impossibility:2019,rohde.zeiler:riemann:2018}.
This can be achieved by imposing further algebraic relations that prescribe the entropy dissipation rate across the interface \cite{truskinovsky:kinks:1993}.
However, determining the explicit form of this relation remains a largely unsolved issue for complex flow regimes, including in particular multi-component flow.
To overcome this fundamental modelling problem on the continuum scale and to avoid using ad-hoc closures, one can look at smaller spatial scales using a molecular-dynamics (MD) approach.
On this ab-initio level the modelling of a multi-component fluid becomes accessible since each molecule of a component has the same physical properties regardless of the fluid phase it is part of (see Section \ref{ch:md_mix}).
However, MD simulations can only be run for very limited spatial and temporal scales.
We utilize a tailored MD algorithm to solve MD Riemann problems across the (discrete) interface only, whereas the dynamics in the bulk phase domains is approximated by a standard finite-volume method for the first-order two-component system \eqref{ch:multicomponent}.
With this hybrid ansatz we propose a heterogeneous multiscale model (HMM, in the sense of \cite{e.engquist.ea:heterogeneous:2007}) that accounts for physically accurate interfacial mass and momentum transfer, as well as the motion of the phase boundary by an MD interface solver, see Section \ref{ch:multiscale_multicomponent}.
A critical issue in this context is the coupling of the different models, i.e., from MD systems with pairwise interactions to continuum-scale Euler-type equations and vice versa.
We rely on classical, statistical averaging to infer continuum-scale quantities from microscale particle simulations.
In that way, we transfer the microscale behavior at fluid interfaces to the continuum scale via a data-based approach.
Concerning the numerical discretization on the continuum scale we employ a moving-mesh finite-volume method which enables us to include the MD dynamics directly at the phase boundary, while resolving the sharp interface directly within the mesh.
The corresponding moving mesh-algorithms is presented in \cite{chalons.rohde.ea:finite:2017, MagieraRohde22}, see \cite{alkamper2021interface} for the corresponding open-source code.\\
Finally, to overcome the computational costs of expensive MD simulations, we apply a surrogate model for the MD interface solver.
This surrogate model is based on the constraint-aware neural networks that have been developed in \cite{magiera.ray.ea:constraint:2020}, which come up with the same computational efficiency, but ensure mass conservation of the numerical discretization, in particular across the interface.\\
In Section \ref{ch:results:mix}, we present numerical results for the HMM for compressible two-component flow with phase transitions.
We focus on a mixture of the noble gases argon and methane using their exact physical properties.
In particular this includes a specific equation of state that is used in the bulk domains.
The results comprise simulations in one, two, and three space-dimensions and settle, for multiple space-dimensions, around the impingement of pressure waves on droplets inducing oscillation and condensation processes.
We summarize our findings in the concluding Section \ref{sec:summary}.

Up to our knowledge, the continuum-scale phase-transition regimes in Section \ref{ch:results:mix} for realistic argon--methane mixtures have not been computed before using a physically accurate model for the interfacial motion.
This contribution is the main novelty of the paper at hand, showing that the approach we have developed for the much simpler single-component cases in \cite{magiera.rohde:particle:2018,magiera.rohde:analysis:2021,MagieraRohde22} can be readily extended to complex scenarios lacking any knowledge about the proper choice of transmission conditions on the continuum scale.

In the remainder of the introductory section we embed our work into the existing literature for sharp-interface modelling in phase-transition flow.
In fact, two-phase fluid flow models with sharp interfaces for a single component are widely investigated, see e.g. \cite{fechter.munz.ea:sharp:2017,saurel.petitpas.ea:modelling:2008,faccanoni.kokh.ea:modelling:2012,ma.lv.ea:entropy:2017,ghazi.james.ea:nonisothermal:2021,shen.ren.ea:3d:2020} for various modelling approaches in the mathematical realm.
As mentioned above, much less is known for multi-component flow.
We refer to the works \cite{HMR:16,hantke.mueller:analysis:2018} which address the closure
problem but do not provide a full continuum-mechanical closure.
An alternative ansatz \cite{HMWY:23} builds on the integration of first-order models into phase-field models \cite{DreyerGiesselmann:14}.\\
MD-simulations are readily and since long time used tools for the dynamics of phase boundaries.
However, multiscale modelling that connects to the continuum-scale dynamics for compressible liquid--vapor flow has not been done up to our knowledge.
The works in e.g. \cite{frezzotti.barbante:simulation:2020,hitz.joens.ea:comparison:2020} indeed support the excellent match of up-scaled quantities if an appropriate closed continuum-scale model exists.
\section{Modelling of Two-Component, Two-Phase Flow}
\label{ch:multicomponent}
In this section, we first review a continuum-scale model from \cite{bothe.dreyer:continuum:2015} that governs the description of two-phase flow for a fluid that consists of several components.
We focus on argon--methane mixtures with comparably well-known fluid properties.
Albeit the fact that the governing equations take the form of first-order balance laws, it turns out that the analytical understanding of hyperbolicity and the properties of solutions for e.g.\ Riemann problems is quite limited, and cannot serve as the basis for numerical simulation.\\
Throughout this work, we consider all variables in combination with their physical units.
This is important for converting microscale quantities to continuum-scale quantities later on.
Moreover, the simulations in Section \ref{ch:results:mix} refer to realistic fluid regimes for the argon--methane mixture.

\subsection{An Isothermal Two-Component Flow Model}
\label{sec:iso_mix_model}
We consider isothermal, inviscid, non-reactive fluids consisting of two components, which we mark by an index \( a\in \{\cmpa,\allowbreak \cmpb\}\).
Later the index $\cmpa$ indicates argon and the index $\cmpb$ methane.
\newline
To model two-component flow, we choose the so-called \mbox{class-II} model for multi-component flow which has been derived in \cite{bothe.dreyer:continuum:2015} as a thermodynamically consistent ansatz.
It includes Maxwell--Stefan diffusion to account for the crucial frictional interaction between components.
The model can be written in terms of a first-order balance law.
Precisely, we have for each component \(a\in \{\cmpa,\allowbreak \cmpb\}\) the equations
\begin{align}
    \label{eq:iso_euler_multicomponent}
    \begin{aligned}
        \partial_t \rho_a + \nabla \cdot (\rho_a \vec{v}_a)                                                     & = 0,                                                                \\
        \partial_t (\rho_a \vec{v}_a) + \nabla \cdot (\rho_a \vec{v}_a \otimes \vec{v}_a) + \rho_a \nabla \mu_a & = - T \frictionvariable_{ab} \rho_a \rho_b (\vec{v}_a - \vec{v}_b),
    \end{aligned}
\end{align}
in \(\contdomain_\phaserl(t)\) for \(t \in (0,t_{\mathrm{end}})\).
The primary variables are the partial mass densities \(\rho_a\) [\si{\kilogram\per\cubic\metre}] and the partial velocities \(\vec{v}_a\) [\si{\metre\per\second}].
The constant reference temperature of the mixture is denoted by \(T\) [\si{\kelvin}].
The function \(\mu_a = \mu_a(\rho_{\cmpa}, \rho_{\cmpb}; T)\) [\si{\joule\per\kilogram}] represents the chemical potential of component \(a\in \{\cmpa,\allowbreak \cmpb\}\) with respect to the partial mass densities and temperature.
We will discuss a specific choice for $\mu_a$ below.
\newline
The term \(\frictionvariable_{\cmpa\cmpb} = \frictionvariable_{\cmpb\cmpa} > 0\) denotes the friction factor between the components \(\cmpa\), \(\cmpb\) and is proportional to the reciprocal of the Maxwell--Stefan diffusion coefficients \(\mathDJ_{\cmpa\cmpb}\) [\si{\metre\squared\per\second}].
More specifically, \(\frictionvariable_{\cmpa\cmpb}\) computes as
\begin{align}
    \frictionvariable_{\cmpa\cmpb} = \frac{\mathrm{R}}{M_{\cmpa} M_{\cmpb} \, c} \cdot \frac{1}{\mathDJ_{\cmpa\cmpb}}.
\end{align}
Here, \(c\) [\si{\mole\per\cubic\metre}] denotes the total molar concentration, which is defined by \(c = c_{\cmpa} + c_{\cmpb}\), where \(c_{a} = \rho_a/M_a\) [\si{\mole\per\cubic\metre}] are the molar concentrations, and \(M_a\) [\si{\kilogram\per\mole}] the molar masses of component \(a\in \{\cmpa,\allowbreak \cmpb\}\).
The ideal gas constant is given by \(\mathrm{R} \approx \SI[allow-number-unit-breaks=false]{8.314}{\joule\per\mole\per\kelvin}\).
Another relevant quantity in the context of two-component flow are the mole fractions \(x_a\) [--], which are defined as
\(x_a = c_a/c\), \(a\in \{\cmpa,\allowbreak \cmpb\}\).
They describe the ratio of the number of \(\cmpa\)/\(\cmpb\)-component particles to the total amount of particles.
For an in-depth discussion of friction in multi-component systems we refer to e.g. \cite{krishna.wesselingh:maxwell:1997}.\\
The initial data of \eqref{eq:iso_euler_multicomponent} are comprised of a connected initial phase boundary \(\continterface(0)\), initial partial mass densities \(\rho_{a,0}(\vx)\), and initial partial velocity fields \(\vv_{a,0}(\vx)\), such that
\begin{align}
    \label{eq:iso_euler_mix_iv}
    \rho_a(\vx,0)    & = \rho_{a,0}(\vx),
                     &
    \vec{v}_a(\vx,0) & = \vv_{a,0}(\vx),
                     & \text{ for } a \in \{\cmpa,\cmpb\} \text{ and } \vx \in \contdomain_{\phaserl}(0).
\end{align}
We define $\UU_0(\vx)=(\rho_{\cmpa,0}(\vx),\rho_{\cmpb,0}(\vx), (\rho_{\cmpa,0}\vv_{\cmpa,0})(\vx)^{\top},(\rho_{\cmpb,0}\vv_{\cmpb,0})(\vx)^{\top})^{\top}$.
Boundary conditions will be specified in Section \ref{ch:results:mix}.

\begin{figure}[htbp]
    \centering
    \includegraphics[width=0.66\columnwidth]{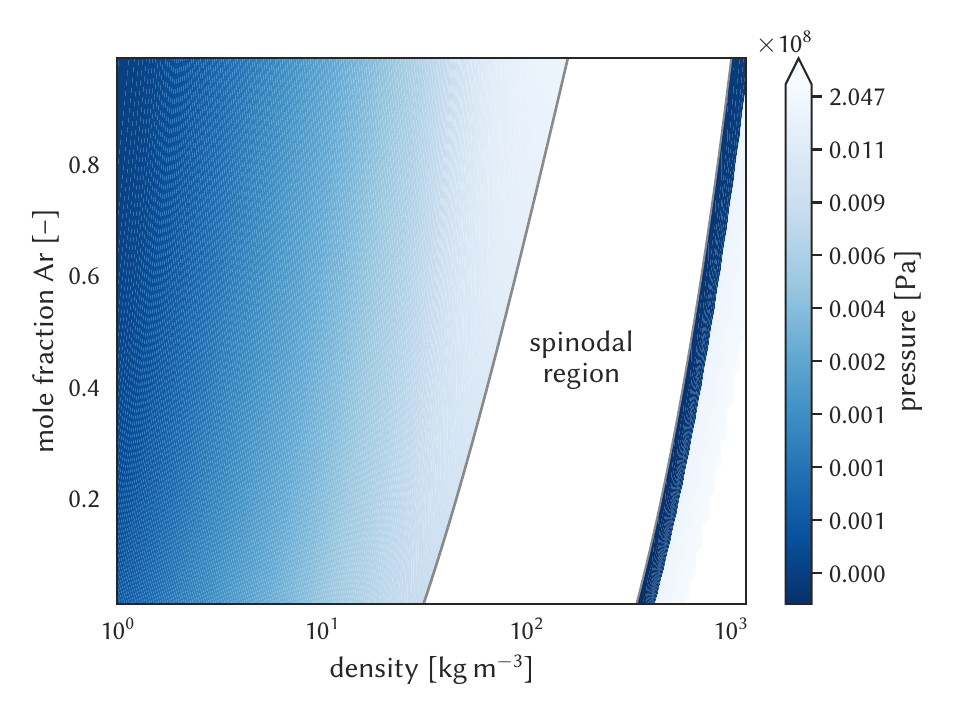}
    \caption{Pressure of the {\pcsaft} {EOS} for an argon--methane mixture at \(T=\SI{110}{\kelvin}\).
        The left-hand (right-hand) region displays vapor (liquid) states.
    }
    \label{fig:eos_pcsaft_ar_me_110}
\end{figure}

To close the system \eqref{eq:iso_euler_multicomponent} we have to provide expressions for the chemical potentials $\mu_{\cmpa},\mu_{\cmpb}$.
This means in the isothermal case that an equation of state (EOS) has to be defined which describes the functional dependencies for the chemical potentials \(\mu_a =
\mu_a(\rho_{\cmpa}, \rho_{\cmpb};T)\).
For two-component flow we choose the {\pcsaft}
    {EOS} \cite{gross.sadowski:perturbed:2001}.
Its underlying theoretical framework provides us in particular with precise parameter values describing argon--methane mixtures, see Section \ref{sec:parameter_table:mix}.
This is the first reason for choosing this mixture in Section \ref{ch:results:mix} on realistic numerical simulations (see Section \ref{ch:md_mix} for the second reason related to molecular-dynamical modelling).
We do not give the highly complex functional form of the {\pcsaft} EOS.
Rather, in Figure~\ref{fig:eos_pcsaft_ar_me_110} we display a plot of the mixture pressure for the argon--methane mixture at $T = \SI{110}{\kelvin}$ as a function of the total density $\rho_0 +\rho_1$ and the
mole fraction $x_{\cmpa}$ of argon.
We observe a vapor-state domain and, separated by a spinodal region, a liquid-state domain.
Although no proof is available we claim that for sufficiently small temperatures the set
\( \{ (\rho_0, \rho_1, \vm_0^{\top}\coloneqq \rho_0\vv_0^{\top}, \vm_1^{\top}\coloneqq \rho_1\vv_1^{\top} )^{\top} \setsep \rho_0,\rho_1>0,\, \vv_0,\vv_1\in \bR^d \}\)
is separated into a liquid-state region ${\mathcal P}_{\phasel}$ and a vapor-state region ${\mathcal P}_{\phasel}$ which define the phases in our case.
Having defined a liquid and a vapor phase, we assume that the two-component fluid in $\contdomain$ partitions the domain in two subdomains $\contdomain_\pm= \contdomain_\pm(t)$ separated by a single, sharp phase boundary \(\continterface(t)\), see Figure~\ref{fig:geometry}.
Just as for the topology of the state space and the splitting, there is up to now no proof that the left-hand side of system \eqref{eq:iso_euler_multicomponent} is hyperbolic in the state space
\[
    {\mathcal U} = {\mathcal P}_{\phasel} \cup {\mathcal P}_{\phaser},
\]
and elliptic in the spinodal region.
Motivated by the situation for single-component, two-phase flow this is exactly what we assume from now on \cite{rohde.zeiler:riemann:2018}.
In all the numerical simulations we have not encountered contrary evidence.

\subsection{Planar Two-Phase Solutions and the Riemann Problem}
The system \eqref{eq:iso_euler_multicomponent} consists of first-order evolution laws.
To gain insight into the two-phase behavior of the solution, we fix \(t\in [0,t_{\mathrm{end}}]\) and some point
$\boldsymbol{\xi}$ on the interface $\boldsymbol{\Gamma}(t)$.
We consider planar solutions of \eqref{eq:iso_euler_multicomponent} on the line depicted by the normal \(\vn = \vn(\boldsymbol{\xi},t)\in \bS^{d-1}\) (pointing in the direction of the vapor-state domain $\contdomain_{\phaser}(t)$), see again Figure~\ref{fig:geometry}.\\
We define, for some state \(\UU \in \mathcal U \subset \bR^{2+2d}\), the normal velocities $ v_0 = \vv_0 \cdot \vn$, $v_1 = \vv_1 \cdot \vn$ and the $\vn$-rotated state
\begin{equation}\label{def:projected_state_variable}
    \uu \coloneqq \UU_{\parallel\vn} = (\rho_0,\rho_1, m_0\coloneqq \rho_0v_0,m_1\coloneqq \rho_1 v_1)^{\top} \coloneqq (\rho_0,\rho_1,\rho_0\vec{v}_0\cdot \vn,\rho_1\vec{v}_1\cdot \vn)^{\top} \in \bR^{4}.
\end{equation}
Furthermore, we set
\[
    \vec{q}(\uu) \coloneqq - T \frictionvariable_{01} \rho_0\rho_1 \big(0,0,v_0-v_1, v_1-v_0 \big)^\top.
\]
Then, the rotational invariance of the Euler-type system \eqref{eq:iso_euler_multicomponent} ensures that there is a matrix $\mat{A} = \mat{A}(\uu;\vn) \in \bR^{4\times 4} $ such that the solution $\uu=\uu(x,t)$ of the $\vn$-rotated (one-dimensional) problem
\begin{align} \label{eq:rotated_cl}
    \begin{aligned}
        \partial_t \uu + \mat{A}(\uu;\vn) \partial_x \uu & = \vec{q}(\uu),
                                                         &                                          & \text{ in } \bR \times (0,t_{\mathrm{end}}), \\
        \uu(\cdot,0)                                     & = \uu_0 \coloneqq \UU_{0,\parallel\vn} ,
                                                         &                                          & \text{ in } \bR,
    \end{aligned}
\end{align}
defines a planar solution of \eqref{eq:iso_euler_multicomponent}.
We recover the original states in $\mathcal U$ by the operator
\begin{align} \label{def:projected_state_variables}
    P_\vn \colon (\rho_0,\rho_1,m_0,m_1)^{\top} \mapsto (\rho_0,\rho_1, m_{0,1} n_1,\ldots, m_{0,d} n_d, m_{1,1} n_1,\ldots, m_{1,d}n_d,)^{\top},
\end{align}
with $n_1,\ldots,n_d$ being the entries of $\vn$.
Now, posing the Riemann problem for the system \eqref{eq:rotated_cl} means to choose initial data of the form
\begin{align} \label{Riemannini}
    \uu_0(x) & = \begin{cases}
                     \uu_{\phasel} & \text{ for } x < 0, \\
                     \uu_{\phaser} & \text{ for } x > 0,
                 \end{cases}
\end{align}
where $\uu_\pm$ are computed from given Riemann states \(\UU_{\pm} \in \mathcal U\) by
\[
    \uu_\pm \coloneqq (\UU_{\phaserl})_{\parallel\vn}.
\]
The Riemann states \(\UU_{\pm} \in {\mathcal P}_\pm \) should be understood as traces states in $\contdomain_\pm(t)$ respectively.\\
The Riemann problem for the $\vn$-rotated system \eqref{eq:rotated_cl} will play a crucial role in the construction of the heterogeneous multiscale method (HMM) in Section~\ref{ch:multiscale_multicomponent}.
The most convenient situation would be to have an exact Riemann solver for \eqref{eq:rotated_cl}\textsubscript{1}.
Motivated by the one-phase case one can then expect that the solution of the Riemann problem \eqref{eq:rotated_cl}, \eqref{Riemannini} includes exactly one isolated phase boundary connecting two constant states in ${\mathcal P}_\pm$.
Its speed $s\in\bR$ and the values of the adjacent states
$ \uu^*_{\pm} \coloneqq (\rho^{*}_{\cmpa,\pm}, \rho^{*}_{\cmpb,\pm},
    v^{*}_{\cmpa,\pm}, v^{*}_{\cmpb,\pm})^{\top} \in {\mathcal P}_\pm $
would then provide all the necessary information about the future behavior of the solution locally at $ \vec{\xi}\in \boldsymbol{\Gamma}(t)$.
Let us denote such a (perfect) interface solver in terms of the original variables for \eqref{eq:iso_euler_multicomponent}, i.e., using the mapping \eqref{def:projected_state_variable} and the back projection $ P_\vn$, by
\begin{align}
    \label{eq:multiscale_mix_riemannsolver}
    \genericsolversym \colon (\UU_{-},\UU_{+};\vn) \mapsto (\UU^{\ast}_{-},\UU^{\ast}_{+}, s).
\end{align}
Actually, an exact Riemann solver $\genericsolversym \colon \PP_{\phasel} \times \PP_{\phaser} \times \bS^{d-1} \to \PP_{\phasel} \times \PP_{\phaser} \times \bR$ for arbitrary two-phase Riemann data $\UU_\pm$
is not available for multicomponent flow (but see \cite{hantke.mueller:analysis:2018,HMWY:23} for some works in this direction).
This problem is of course related to the lacking information about hyperbolicity of \eqref{eq:iso_euler_multicomponent}.
We aim to substitute the missing Riemann solver by an interface solver based on MD simulations which will be discussed in the next section.
\section[MD: Multicomponent Fluids]{Molecular Dynamics for Multicomponent Fluids}
\label{ch:md_mix}
One big advantage of molecular dynamics (MD) are their versatility to simulate fluid mixtures consisting of more than one fluid component.
Fluid components in this context denote different types of molecules that make up the composition of a fluid.
\newline
In this section, we review the molecular-dynamical modelling for binary mixtures.
Then we extend a numerical method that has been developed in \cite{MagieraRohde22} for one-component fluids to the case of binary mixtures, i.e., fluid mixtures consisting of two components.
We also discuss the compatibility of MD modelling for argon--methane mixtures with the \pcsaft{} EOS introduced in Section~\ref{sec:iso_mix_model}.
As before, we indicate the two components by the numbers \(\cmpa\)~and \(\cmpb\).

\subsection{Molecular Dynamics for Binary Mixtures}
\begin{figure}[tb]
    \centering
    \includegraphics[width=0.5\columnwidth]{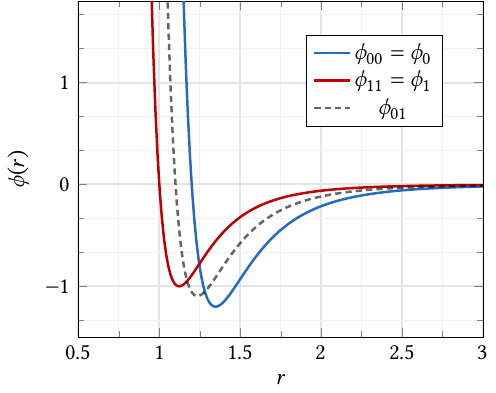}
    \caption{Two-component Lennard--Jones potential for the parameters
        \(\varepsilon_{\cmpa} = 1.2\),
        \(\sigma_{\cmpa} = 1.2\), and
        \(\varepsilon_{\cmpb} = 1\),
        \(\sigma_{\cmpb} = 1\).
        The two solid-line graphs display intra-component interaction, and the dashed one the interaction potential between the two components, using the combination rules \protect\eqref{eq:lorentz-berthelot}.
    }%
    \label{fig:lj_mix_potential}
\end{figure}

We focus on fluid components that can be modeled by spherical, rotational invariant Lennard--Jones particles, such as noble gases or relatively simple molecules like methane.
In the two-component setting, each particle belongs to either component~\cmpa, or to component~\cmpb.
We denote the particle type for the \(i\)-th particle by \(\varsigma_i\).
Accordingly, the mass of each particle is \(m_{\cmpa}\) [\si{\atomicmassunit}], if \(\varsigma_i = \cmpa\), or alternatively \(m_{\cmpb}\) [\si{\atomicmassunit}], if \(\varsigma_i = \cmpb\). \\
On the atomistic scale the dynamics of an ensemble of $N= N_{\cmpa}+N_{\cmpb}$ particles with \(N_{\cmpa} \in \bN\) particles of component \(\cmpa\) and \(N_{\cmpb} \in \bN\) particles of component \(\cmpb\) are governed by the ordinary initial value problem \cite{allen.tildesley:computer:2017}
\begin{align}\label{mdequations}
    \begin{aligned}
        \vx_i'(\tau) & = \vv_i(\tau),
                     &                                                                                                                        %
        \vv_i'(\tau) & = \va_i(\tau),
                     &                                                                                                                        %
        \va_i(\tau)  & = - \sum_{j \neq i} \nabla_{\vx_i} \phi_{\varsigma_i \varsigma_j}(\norm{\vx_i(\tau) - \vx_j(\tau)}) / m_{\varsigma_i},
        \\
        \vx_i(0)     & = \vx^{0}_{i},
                     &
        \vv_i(0)     & = \vv^{0}_{i},
    \end{aligned}
\end{align}
for the particle positions \(\vx_i(\tau) = (x_i(\tau), y_i(\tau), z_i(\tau))^\top \in \bR^3\), velocities \(\vv_i(\tau) \in \bR^3\), and accelerations \(\va_i(\tau) \in \bR^3\) of the \(i\)-th particle, with \(\tau \in [0,\tau_{\mathrm{end}}]\). Here $\tau_{\mathrm{end}}>0$ is the MD end time.
The MD simulations are carried out in the atomistic-scale units [\si{\kelvin}] for energy, [\si{\angstrom}] for distances, and [\si{\atomicmassunit}] for mass.
The crucial ingredients for the MD model \eqref{mdequations} are the Lennard--Jones potentials $\phi_{ab}$ with
\begin{align} \label{eq:lj_mix_potential}
    \phi_{ab}(r) = 4 \varepsilon_{ab} \left( \left( \frac{\sigma_{ab}}{r} \right)^{12} - \left( \frac{\sigma_{ab}}{r} \right)^6 \right).
\end{align}
The Lennard--Jones parameters \(\varepsilon_{\cmpa\cmpb} = \varepsilon_{\cmpb\cmpa}\), \(\sigma_{\cmpa\cmpb} = \sigma_{\cmpb\cmpa}\) depend on the corresponding parameters \(\varepsilon_{\cmpa}\), \(\sigma_{\cmpa}\), \(\varepsilon_{\cmpb}\), \(\sigma_{\cmpb}\) for single-component potentials.
If both interacting particles belong to the same component $a\in \{0,1\}$, the parameters are given by \(\varepsilon_{aa} \coloneqq \varepsilon_{a}, \sigma_{aa} \coloneqq \sigma_{a}\).
If they belong to different components, the Lorentz--Berthelot combination rules \cite{lorentz:ueber:1881,berthelot:sur:1898} suggest for some scaling numbers \(\eta\) [--] and \(\xi\) [--]
\begin{align}
    \label{eq:lorentz-berthelot}
    \sigma_{\cmpa\cmpb} = \sigma_{\cmpb\cmpa}           & = \eta \frac{\sigma_{\cmpa} + \sigma_{\cmpb}}{2},
                                                        &
    \varepsilon_{\cmpa\cmpb} = \varepsilon_{\cmpb\cmpa} & = \xi \sqrt{\varepsilon_{\cmpa} \varepsilon_{\cmpb}}.
\end{align}
In Figure \ref{fig:lj_mix_potential}, the Lennard--Jones potentials for different parameters are plotted.

As discretization for \eqref{mdequations} we use the classical Velocity-Verlet algorithm \cite{hairer.lubich.ea:geometric:2003} which amounts to performing the following steps with a given time step \(\Delta \tau > 0\), \(N_{\mathrm{end}} = \lceil \tau_{\mathrm{end}} / \Delta \tau \rceil\) times:
\begin{align}\label{VValgo}
    \begin{aligned}
        \text{1. } &
        \vx^{n+1}_i = \vx^{n}_i + \Delta \tau \, \vv^{n}_i + \tfrac{1}{2} \Delta \tau^2 \, \va^{n}_i \text{ for all } i = 1, \ldots, N. \\
        \text{2. } &
        \vv^{n+\sfrac{1}{2}}_i = \vv^{n}_i + \tfrac{1}{2} \Delta \tau \, \va^{n}_i \text{ for all } i = 1, \ldots, N.                   \\
        \text{3. } &
        \text{Compute accelerations } \va^{n+1}_i \text{ for all } i = 1, ,\ldots, N \text{ as in \eqref{mdequations}},                 \\ & \text{using the positions } \vx^{n+1}_j. \\
        \text{4. } & \vv^{n+1}_i = \vv^{n+\sfrac{1}{2}}_i + \tfrac{1}{2} \Delta \tau \, \va^{n+1}_i \text{ for all } i = 1, \ldots, N.
    \end{aligned}
\end{align}

\begin{remark}[Reduced units, time scales, and cutoff-potential]\mbox{}\newline
    \begin{itemize}
        \item[(i)] The MD simulations using the Velocity-Verlet algorithm \eqref{VValgo} are done in reduced units according to \autoref{table:reduced_units} below.
            With this procedure we follow \cite{allen.tildesley:computer:2017}. It avoids floating point over-/underflow in performing the numerical algorithm.\\
            The number of particles $N$, the time step $\Delta \tau$ (in reduced units) and the number of time steps $N_{\text{end}}$
            will enter among others as parameters in Algorithm \ref{alg:md_mix_riemann} for solving an MD Riemann problem.
            The precise values of the method's parameters are listed in Table \ref{tab:mix_paramater}.
            These parameters are chosen in a way that
            allows a reliable extraction of Riemann data in the final multiscale method which acts on the continuum scale. Note that we expect the MD Riemann problem
            to have a self-similar solution such that we can extract data at any point in
            time. Without any information about numerical errors
            the ad-hoc choice of these parameters remains a weakness of the entire
            multiscale approach.
        \item[(ii)]
            The application of the algorithm \eqref{VValgo} leads to a quadratic complexity in terms of the total particles number $N$.
            This is due to the particle-wise interaction in \eqref{eq:lj_mix_potential}.
            A remedy are cutoff-potentials.
            For a two-component mixture we substitute $\phi_{ab}$ for \(a,b \in \{\cmpa, \cmpb\}\) by
            \begin{align} \label{eq:mix_cutoff_potential}
                \phi_{ab}(r; r_{\mathrm{cutoff}})
                \coloneqq \begin{cases}
                              \phi_{ab}(r) & : \,\,r < \sigma_{ab} r_{\mathrm{cutoff}},     \\
                              0            & : \,\, r \geq \sigma_{ab} r_{\mathrm{cutoff}}.
                          \end{cases}
            \end{align}
            Like all numerical parameters for the MD simulations the actually used value for the cutoff radius $r_{\mathrm{cutoff}} $ can be found in Table \ref{tab:mix_paramater} in the appendix.
            The error introduced by the cutoff functional will be corrected by the introduction of long-range interaction forces on the discrete level.
            We refer to \cite{magiera:molecular:2021} for a special choice for two-component flows.
    \end{itemize}
\end{remark}
\begin{table}[htb]
    \centering
    \begin{tabular}{lrr}
        \toprule
                    & \textbf{Reduced units}                                                & \textbf{SI units}                                                                  \\
        \midrule
        mass        & \SI{1}{\atomicmassunit}                                               & \SI[round-mode=places,round-precision=2]{1.660539040e-27}{\kilogram}               \\
        length      & \SI{1}{\redlength}                                                    & \SI{e-10}{\metre}                                                                  \\
        energy      & \SI{1}{\redenergy}                                                    & \SI[round-mode=places,round-precision=2]{1.380e-23}{\joule}                        \\
        \midrule
        velocity    & \SI{1}{\redenergy\tothe{0.5}\per\atomicmassunit\tothe{0.5}}           & \SI[round-mode=places,round-precision=2]{91.1622421005}{\metre\per\second}         \\
        time        & \SI{1}{\redlength\atomicmassunit\tothe{0.5}\per\redenergy\tothe{0.5}} & \SI[round-mode=places,round-precision=2]{1.0969454e-12}{\second}                   \\
        pressure    & \SI{1}{\redenergy\per\cubic\redlength}                                & \SI[round-mode=places,round-precision=2]{1.380e+7}{\pascal}                        \\
        temperature & \SI{1}{\redenergy\per\boltzmannunit}                                  & \SI{1}{\kelvin}                                                                    \\
        density     & \SI{1}{\atomicmassunit\per\cubic\redlength}                           & \SI[round-mode=places,round-precision=2]{1.66053904e+3}{\kilogram\per\cubic\metre} \\
        \bottomrule
    \end{tabular}
    \captionsetup{width=0.75\textwidth}
    \caption{Unit conversion table.}
    \label{table:reduced_units}
\end{table}
We are interested in running MD-simulations that provide an approximation of the solution of the
Riemann problem \eqref{eq:rotated_cl}, \eqref{Riemannini}.
For this purpose the $\vn$-rotated Riemann initial data $ \uu_\pm$ from \eqref{Riemannini} have to be represented as initial data for \eqref{mdequations}, i.e., we need an initial particle and velocity distribution in 3D.
In \cite{MagieraRohde22} we have suggested a linked-cell algorithm to realize this transfer of single-component, two-phase constant Riemann data into a cubic domain $\Sigma$ that consists of two cuboids $\Sigma_\pm$ separated by the atomistic interface.
The algorithm includes a thermalization process that makes sure that the MD-particle and
particle velocity distributions matches the constant reference temperature $T$, see \eqref{eq:md_mix:temperature} below for
the definition of atomistic-scale temperatures.
The linked-cell algorithm can be readily applied to the two-component case and is not further detailed.

In the next step, the system \eqref{mdequations} is solved numerically by the scheme \eqref{VValgo}.
We can again utilize a procedure from \cite{MagieraRohde22} to determine the approximate position of the atomistic interface and in particular an approximate interface speed $s\in \bR$ like in the
exact interface solver $\mathcal R$ from \eqref{eq:multiscale_mix_riemannsolver}.
To obtain all output data in \eqref{eq:multiscale_mix_riemannsolver} we need to find approximations of adjacent continuum-scale states $\UU^*_\pm$, $\uu^*_\pm$, respectively.
Thus, we need to define continuum-scale quantities from the MD-simulation data.
\newline
Let us consider a system of \(N\) particles.
The component-wise mass concentrations \(\widehat{\rho}_{\cmpa}(\vx,\mdtimevar)\), \(\widehat{\rho}_{\cmpb}(\vx,\mdtimevar)\)
[\si{\atomicmassunit\per\cubic\angstrom}],
and
momentum concentrations \(\widehat{\vm}_{\cmpa}(\vx,\mdtimevar)\),
\(\widehat{\vm}_{\cmpb}(\vx,\mdtimevar)\)
[\si{\kelvin\tothe{0.5}\atomicmassunit\tothe{0.5}\per\cubic\angstrom}] can be computed by the Irving--Kirkwood formulas in the following way, for \(a \in \{\cmpa, \cmpb\}\),
\begin{align}\label{IK1}
    \widehat{\rho}_{a}(\vx,\mdtimevar) = \sum_{\substack{i=1 \\ \varsigma_i = a}}^{N} m_i \, \delta_{\vx_i(\mdtimevar)}(\vx),\qquad
    \widehat{\vm}_{a}(\vx,\mdtimevar) = \sum_{\substack{i=1  \\ \varsigma_i = a}}^{N} m_i \, \vv_i(\mdtimevar) \, \delta_{\vx_i(\mdtimevar)}(\vx),
\end{align}
where \(m_i\), \(\vx_i\), \(\vv_i\) are the mass, position and velocity of the \(i\)-th particle.
The local temperature distribution of the mixture can be formally understood as
\begin{align}
    \label{eq:md_mix:temperature}
    \widehat{T}(\vx,\mdtimevar) = \frac{1}{3} \left( \sum_{i=1}^{N} \delta_{\vx_i(\mdtimevar)}(\vx) \right)^{-1} \sum^{N}_{i = 1} (m_i (\vv_i - \overline{\vv}) \cdot (\vv_i - \overline{\vv})) \, \delta_{\vx_i(\mdtimevar)}(\vx),
\end{align}
where \(\overline{\vv} = \overline{\vv}(\vx,\mdtimevar)\) [\si{\kelvin\tothe{0.5}\atomicmassunit\tothe{-0.5}}] denotes the local barycentric average velocity of the mixture.
The local average velocity can be realized by e.g.
\begin{align}
    \label{eq:continuous_avg_velocity:md}
    \overline{\vv}(\vx,\mdtimevar) =
    \Bigl( \sum_{i=1}^{N} K^{\varepsilon}(\abs{\vx - \vx_i(\mdtimevar)}) \Bigr)^{-1}
    \Bigl( \sum_{i=1}^{N} \vv_i(\mdtimevar) K^{\varepsilon}(\abs{\vx - \vx_i(\mdtimevar)}) \Bigr),
\end{align}
where \(K^{\varepsilon} \colon \bR \to \bR \) is a smooth kernel with length-scale parameter \(\varepsilon > 0\).
It cannot be expected that the temperature field $ \widehat{T}$ from \eqref{eq:md_mix:temperature} equals the constant temperature $T$ in the isothermal system \eqref{eq:iso_euler_multicomponent}.
This can be partially circumvented by applying a thermostat, but we will observe in Section~\ref{ch:results:mix} that substantial discrepancies can occur for the vapor-phase region.\\
Using the formulas \eqref{IK1} we obtain the adjacent states by sampling over small spatial slices
$\Sigma^*_\pm$ next to the interface.
That means we compute for \(a \in \{\cmpa, \cmpb\}\)
\begin{align} \label{eq:md_mix_interface_states}
    \begin{aligned}
        \rho^{*}_{a,\phaserl}(\mdtimevar)
         & = \operatorname{vol}(\Sigma^*_{\phaserl})^{-1} \sum_{{x_i \in \Sigma^*_{\phaserl}, \, \varsigma_i = a}} m_i,
        \\
        \vv^{*}_{a,\phaserl}(\mdtimevar)
         & = \abs{x_i \in \Sigma^*_{\phaserl}, \varsigma_i = a}^{-1} \sum_{{x_i \in \Sigma^*_{\phaserl}, \, \varsigma_i = a}} \vv_i,
    \end{aligned}
\end{align}
where \(\operatorname{vol}(\Sigma^*_{\phaserl}) \) denotes the volume of the sampling region \(\Sigma^*_{\phaserl}\), and \(\abs{x_i \in \Sigma^*_{\phaserl}, \varsigma_i = a}\) the number of component-\(a\) particles inside it.
We can obtain the local interface states for every time step during an MD simulation.
For the final outcome we perform time-averaging of these instantaneous states over a fraction \(\fractionvar_{\mdtimevar\mathrm{-smpl}}\) of the total simulation time \(\mdtimeend\). Note that the
conversion of MD quantities and continuum-scale quantities has to account for the unit conversion given in Table \ref{table:reduced_units}.

We summarize all steps for solving an MD-Riemann problem in the following algorithm,
see also Algorithm~\ref{alg:riemann_mix_euler} in the context of the multiscale method.

\begin{boxalgorithm}[label=alg:md_mix_riemann]{Atomistic-Scale Interface Solver for Binary Mixtures}
    \textbf{Input:}
    initial states \((\rho_{\cmpa,\phasel}\), \(\rho_{\cmpb,\phasel}\),
    \(\vm_{\cmpa,\phasel}\), \(\vm_{\cmpb,\phasel}) \in \PP_{\phasel}\), \((\rho_{\cmpa,\phaser}\), \(\rho_{\cmpb,\phaser}\),
    \(\vm_{\cmpa,\phaser}\), \(\vm_{\cmpb,\phaser}) \in \PP_{\phaser}\),
    normal direction \(\vn \in \bS^{d-1}\).

    \textbf{Parameters:}
    total number of particles \(N_{\mathrm{particles}}\),
    number of time steps \(\mdstepend\),
    time step \(\Delta \mdtimevar > 0\),
    time sampling ratio \(\fractionvar_{\mdtimevar\mathrm{-smpl}} \in [0,1]\),
    sampling regions $\Sigma^*_\pm$.
    \tcbsubtitle{Algorithm}
    \begin{itemize}
        \item Compute initial particle configuration with linked-cell algorithm from \cite{MagieraRohde22}
              starting from the $\vn$-rotated Riemann data $\uu_{\pm} = (\rho_{\cmpa,\pm}, \rho_{\cmpb,\pm},
                  v_{\cmpa,\pm}, v_{\cmpb,\pm})^{\top} $.
        \item For \(n = 1,\ldots, \mdstepend\):
              \begin{itemize}
                  \item Run method \eqref{VValgo}, using the interaction potential \eqref{eq:mix_cutoff_potential}.
                  \item Compute the interface positions \(\mdinterface(\mdtimevar_{n})\) and interface speeds $s(\mdtimevar_{n})$.
                  \item Compute continuum-scale states \((\rho^{*}_{0,\phaserl}, \rho^{*}_{1,\phaserl}, \vv^{*}_{0,\phaserl}, \vv^{*}_{1,\phaserl}) (\mdtimevar_{n})\) like in \eqref{eq:md_mix_interface_states}.
              \end{itemize}
        \item Compute the time averaged values
              \((\rho^{*}_{\cmpa,\phaserl}, \allowbreak \rho^{*}_{\cmpb,\phaserl}, \allowbreak v^{*}_{\cmpa,\phaserl}, \allowbreak v^{*}_{\cmpb,\phaserl})\),
              \(s\)
              from
              \((\rho^{*}_{\cmpa,\phaserl}, \allowbreak \rho^{*}_{\cmpb,\phaserl}, \allowbreak
              \vv^{*}_{\cmpa,\phaserl}, \allowbreak \vv^{*}_{\cmpb,\phaserl})(\mdtimevar_{n})\), \allowbreak
              \(s(\mdtimevar_{n})\)
              for all \(n = \widetilde{\fractionvar}_{\mdtimevar\mathrm{-smpl}} \mdstepend, \ldots, \mdstepend\).
    \end{itemize}
    \textbf{Result:} interface states \((\rho^{*}_{\cmpa,\phasel}, \allowbreak \rho^{*}_{\cmpb,\phasel}, \allowbreak v^{*}_{\cmpa,\phasel} \allowbreak, v^{*}_{\cmpb,\phasel})\),
    \((\rho^{*}_{\cmpa,\phaser}, \allowbreak \rho^{*}_{\cmpb,\phaser}, \allowbreak v^{*}_{\cmpa,\phaser}, \allowbreak v^{*}_{\cmpb,\phaser})\),
    and interface speed \(s\).
\end{boxalgorithm}
Note that the states $v^*_{a,\pm}(\tau_n)$ are computed from $ \vv^*_{a,\pm}(\tau_n)$ by the back projection $P_\vn$, and that $\tau_n= n\Delta \tau$, $n\in \{0,\ldots, N_{\mathrm{end}}\}$.
\newpage
\subsection{Notes on the Argon--Methane Mixture}
\label{sec:argon-methane}
Up to now, we considered a generic binary mixture.
In the following, we focus on two-component mixtures consisting of argon and methane.
In this case, parameters for the Lennard--Jones potential \eqref{eq:lj_mix_potential} can be found in \cite{VrabecFischer95}, and even experimental data are available \cite{ponte.streett.ea:experimental:1981}.
For the Lennard--Jones parameters in \eqref{eq:lorentz-berthelot} we take the values from \cite{VrabecFischer95} listed in \autoref{tab:md_ar_me_parameters}.
The corresponding combination parameters in \eqref{eq:lorentz-berthelot} are \(\eta = \num{1.00141}\) and \(\xi = \num{0.96400}\), see also
\cite{VrabecFischer95}. For methodological background on the used
approach we refer to \cite{vrabec.lotfi.ea:vapour:1995}.

\begin{remark}[Temperature parameterization of Lennard--Jones parameters]
    The chosen Lennard--Jones parameters are given in \cite{VrabecFischer95} for temperatures $ T= \SI{111}{\kelvin}$ (Argon) and $T= \SI{141}{\kelvin}$ (Methane). The combination parameters refer to the
    temperature $T=\SI{115}{\kelvin}$ and have also been computed in \cite{VrabecFischer95}. In this work the quantities have been shown to be in good fit with available experimental data in \cite{ponte.streett.ea:experimental:1981}. Up to our knowledge
    this setting is commonly used for temperatures well below the critical temperature of the mixture. Notably, our HMM method applies only for settings where the phase boundary can be expected to be a sharp interface separating the liquid and the vapor regions. Close to the critical temperature this separation
    property would fail.
\end{remark}
The studies in the literature refer solely to the homogeneous case, i.e., they do not span the
entire two-phase state space.
In \cite{magiera:molecular:2021} it has been shown by extended MD simulations for homogeneous data covering the whole state space, that the predictions of the {\pcsaft} theory for the chemical potentials $\mu_0,\mu_1$ are in good quantitative
agreement with the MD simulations using the Lennard--Jones potentials for
argon--methane mixtures.
This justifies to use the \pcsaft{} chemical potentials for the continuum-scale computations in the HMM, that will be proposed in Section \ref{ch:multiscale_multicomponent}.
\begin{table}[htb]
    \centering
    \begin{tabular}{
            cc
            S[round-mode=places, round-precision=3, table-format=1.3]
            S[round-mode=places, round-precision=2, table-format=3.2]
            S[round-mode=places, round-precision=3, table-format=2.3]
        }
        \toprule
        \multicolumn{2}{c}{Component} & \(\sigma\) [\si{\angstrom}] & \(\varepsilon / k_{\mathrm{B}}\) [\si{\kelvin}] &
        \(m\) [\si{\atomicmassunit}]
        \\ \midrule
        Argon                         & Ar                          & 3.3967                                          & 117.05 & 39.948 \\
        Methane                       & CH\textsubscript{4}         & 3.7275                                          & 148.99 & 16.043 \\
        \bottomrule                                                                                                                     \\
    \end{tabular}
    \captionsetup{width=\textwidth}
    \caption{Lennard--Jones parameters from \cite{VrabecFischer95} and species masses in reduced units for the argon--methane mixture, cf.\ Table \ref{table:reduced_units}.\label{tab:md_ar_me_parameters}}
\end{table}

\section{A Heterogeneous Multiscale Method for Isothermal Two-Component, Two-Phase Flow}
\label{ch:multiscale_multicomponent}
Following \cite{MagieraRohde22} for single-component flow with phase transitions we introduce a corresponding heterogeneous multiscale method (HMM) that combines the isothermal two-component flow model on the continuum scale (\ref{ch:multicomponent}) with two-component MD simulations described in Section~\ref{ch:md_mix} for the microscale description of the interface motion.
\newline
Here, we focus again on a fluid mixture consisting of argon and methane, see Section~\ref{ch:multicomponent}.
The component-specific quantities are marked by the indices \(\cmpa\) for argon, and \(\cmpb\) for methane.
The two fluid phases are again denoted by the subscripts \(\phasel\) and \(\phaser\) for the liquid and vapor phase.
\newline
First, we formulate on the basis of Algorithm~\ref{alg:md_mix_riemann} the MD-based microscale interface solver $\rmdmix$ which provides an approximate solution for Riemann initial data (Section~\ref{sec:mix_microscale_solver}).
Using an MD-based interface solver directly in place is still too computational expensive.
Therefore we suggest the machine-learned surrogate solver $\networkF_{\nnparams} $ in Section~\ref{sec:surrogate_details_mix}, which exploits the fact that the interface solver $\rmdmix$ can be understood as a nonlinear function mapping finite-dimensional Riemann data to a finite-dimensional space including the phase boundary speed
and the state values adjacent to the phase boundary.
Finally, in Section~\ref{sec:HMMfinal} we present the complete HMM for two-component, two-phase flow.
For the sake of simplicity we ignore in the entire section boundary conditions on $\partial \contdomain$.
\subsection{The Atomistic-Scale Interface Solver {$\rmdmix$}}
\label{sec:mix_microscale_solver}
To reduce the input space of the atomistic-scale interface solver from Algorithm~\ref{alg:md_mix_riemann} we exploit that the continuum-scale system \eqref{eq:iso_euler_multicomponent} is invariant with respect to velocity shifts.
We define the barycentric velocities
\begin{align}
    \label{eq:barycentric_reference_velocity}
    \vv_{\phaserl} \coloneqq
    \big({ \vm_{\cmpa,\phaserl} + \vm_{\cmpb,\phaserl} }\big)/\big({ \rho_{\cmpa,\phaserl} + \rho_{\cmpb,\phaserl}}\big),
\end{align}
and choose $\bar v \coloneqq \vv_{\phasel}$ as the reference velocity.
Furthermore, we introduce the relative velocity for each phase
\begin{align}
    \label{eq:relative_mix_velocity}
    \vv_{\mathrm{rel},\phaserl} \coloneqq \vv_{\cmpa,\phaserl} - \vv_{\cmpb,\phaserl}.
\end{align}
Assuming we know the densities \(\rho_{\cmpa,\phaserl}\), \(\rho_{\cmpb,\phaserl}\), we can compute \((\vv_{\cmpa,\phaserl}, \vv_{\cmpb,\phaserl})\) from \((\vv_{\phaserl}, \vv_{\mathrm{rel},\phaserl})\) and vice versa.
In the sequel this conversion is done implicitly to simplify the notation.
\newline
The following algorithm implements the final MD-interface solver
$\rmdmix \colon \PP_{\phasel} \times \PP_{\phaser} \times \bS^{d-1} \to \PP_{\phasel} \times \PP_{\phaser} \times \bR$ substituting the exact Riemann solver $\cR$.
\begin{boxalgorithm}[label= alg:riemann_mix_euler]{MD-Interface Solver $\rmdmix$: Two-Component, Two-Phase Flow}
    \textbf{Input:}
    initial states \((\rho_{\cmpa,\phasel}\), \(\rho_{\cmpb,\phasel}\),
    \(\vm_{\cmpa,\phasel}\), \(\vm_{\cmpb,\phasel}) \in \PP_{\phasel}\), \((\rho_{\cmpa,\phaser}\), \(\rho_{\cmpb,\phaser}\),
    \(\vm_{\cmpa,\phaser}\), \(\vm_{\cmpb,\phaser}) \in \PP_{\phaser}\),
    normal direction \(\vn \in \bS^{d-1}\).

    \tcbsubtitle{Algorithm}
    \begin{enumerate}
        \item Compute the barycentric and relative phase velocities
              \(\vv_{\phaserl}\), \(\vv_{\mathrm{rel},\phaserl}\) according to
              \eqref{eq:barycentric_reference_velocity}, \eqref{eq:relative_mix_velocity}.
        \item Map \((\rho_{\cmpa,\phaserl}, \rho_{\cmpb,\phaserl},
              \vv_{\phaserl}, \vv_{\mathrm{rel},\phaserl})\) in direction of \(\vn\) according to Definition \ref{def:projected_state_variables}, i.e.,
              \begin{align}
                  \label{eq:mix:multiscale:projected_state_variables}
                  \uu_{\phaserl} =
                  (\rho_{\cmpa,\phaserl}, \rho_{\cmpb,\phaserl},
                  v_{\phaserl}, v_{\mathrm{rel},\phaserl})
                   & \coloneqq
                  \UU_{\parallel\vn, \phaserl}
                  = (\rho_\phaserl, \vv_{\phaserl} \cdot \vn, \vv_{\mathrm{rel},\phaserl} \cdot \vn),
              \end{align}
              and define
              \begin{align*}
                  \UU_{\perp\vn, \phaserl}
                   & \coloneqq
                  (0, 0, \vv_\phaserl - (\vv_\phaserl \cdot \vn) \vn, \vv_{\mathrm{rel},\phaserl} - (\vv_{\mathrm{rel},\phaserl} \cdot \vn) \vn ).
              \end{align*}
        \item
              \label{item:mix_microscale_md_solver:reference_velocity::subtract}
              Consider the barycentric \(\phasel\)-phase velocity as the reference velocity \(\overline{v} \coloneqq v_\phasel\) and map it:
              \begin{align}
                  v_\phaserl
                   & \mapsfrom v_\phaserl - \overline{v}.
              \end{align}
        \item
              \label{item:mix_microscale_md_solver}
              Solve the Riemann problem \eqref{eq:rotated_cl} %
              with Riemann data $(\rho_{\cmpa,\phaserl}, \rho_{\cmpb,\phaserl},
                  \vv_{\phaserl}, \vv_{\mathrm{rel},\phaserl})$ using the MD model.
              This means, we run a two-component MD simulation for the initial continuum states $(\rho_{\cmpa,\phasel}, \rho_{\cmpb,\phasel},{\vec{0}}, \vv_{\mathrm{rel},\phasel})$
              and $(\rho_{\cmpa,\phaser}, \rho_{\cmpb,\phaser}, \vv_{\phaser}, \vv_{\mathrm{rel},\phaser})$ -- see Algorithm~\ref{alg:md_mix_riemann}.
              This yields the wave states
              \((\rho^{*}_{\cmpa,\phasel}, \rho^{*}_{\cmpb,\phasel},
              v^{*}_{\cmpa,\phasel}, v^{*}_{\cmpb,\phasel})\), \((\rho^{*}_{\cmpa,\phaser}, \rho^{*}_{\cmpb,\phaser},
              v^{*}_{\cmpa,\phaser}, v^{*}_{\cmpb,\phaser})\)
              and the wave speed \(s\) in normal direction \(\vn\).
        \item
              \label{item:mix_microscale_md_solver:reference_velocity::add}
              Return the reference velocity \(\overline{v}\) to the output:
              \begin{align*}
                  m^{*}_{\cmpa,\phaserl}
                   & \mapsfrom m^{*}_{\cmpa,\phaserl} + \rho^{*}_{\cmpa,\phaserl} \overline{v}, \\
                  m^{*}_{\cmpb,\phaserl}
                   & \mapsfrom m^{*}_{\cmpb,\phaserl} + \rho^{*}_{\cmpb,\phaserl} \overline{v}, \\
                  s
                   & \mapsfrom s + \overline{v}.
              \end{align*}
        \item Project the directional states back to the full state:
              \begin{align*}
                  \UU^*_{\phaserl}\coloneqq (\rho^{*}_{\cmpa,\phaserl}, \rho^{*}_{\cmpb,\phaserl},
                  \vm^{*}_{\cmpa,\phaserl}, \vm^{*}_{\cmpb,\phaserl})
                   & = (\rho^{*}_{\cmpa,\phaserl}, \rho^{*}_{\cmpb,\phaserl},
                  m^{*}_{\cmpa,\phaserl} \vn, m^{*}_{\cmpb,\phaserl} \vn) + \UU_{\perp\vn, \phaserl}
                  \\
                   & = P_{\vn}((\rho^{*}_{\cmpa,\phaserl}, \rho^{*}_{\cmpb,\phaserl},
                  m^{*}_{\cmpa,\phaserl}, m^{*}_{\cmpb,\phaserl})) + \UU_{\perp\vn, \phaserl},
              \end{align*}
              with \(P_{\vn}\) defined as in \eqref{def:projected_state_variables}.
    \end{enumerate}
    \textbf{Result:} $ \rmdmix(\UU_{\phasel},\UU_{\phaser};\vn) = (\UU^*{_{\phasel}}, \UU^*_{\phaser},s)$
\end{boxalgorithm}

Due to the velocity shift described before we can always achieve that the component $\vm_{0,-}$ can be assumed to vanish when Algorithm \ref{sec:mix_microscale_solver}
is executed within Algorithm \ref{alg:riemann_mix_euler}.
\subsection{The Machine-Learned Interface Solver \(\networkF_{\nnparams}\)}
\label{sec:surrogate_details_mix}
To reduce the computational complexity of Algorithm~\ref{alg:riemann_mix_euler} we employ a machine-learned surrogate interface solver \(\networkF_{\nnparams}\) that approximates \(\rmdmix\), i.e.,
\begin{align}
    \networkF_{\nnparams}{ (\UU_{\phasel},\UU_{\phaser};\vn) }
    )
    \approx
    \rmdmix( (\UU_{\phasel},\UU_{\phaser};\vn)
    \approx {\cR } ( (\UU_{\phasel},\UU_{\phaser};\vn).
\end{align}
For the full multiscale model, the surrogate model \(\networkF_{\nnparams}\) replaces the molecular-scale interface solver \(\rmdmix\) in step~\ref{item:mix_microscale_md_solver} of Algorithm~\ref{alg:riemann_mix_euler}. We note that the implementation of \(\networkF_{\nnparams}\) relies on the reduced-dimensional state spaces in the rotated Riemann problem \eqref{eq:rotated_cl} and using
the invariance with respect to the reference velocity $\vv_{\phasel}$ in \eqref{eq:barycentric_reference_velocity}.
\newline
The final surrogate solver offers large computational gains:
It takes only around
\SI[round-mode=places, round-precision=2]{0.102}{\milli\second}
for a single evaluation.
The MD simulations (Algorithm~\ref{alg:riemann_mix_euler}) on the other hand needs \SIrange{14}{17}{\minute}.
\subsubsection{Data Set Generation}
\label{sec:data_generation:mix}

To prepare the surrogate interface solver, we have to generate a data set \(\datasetvar\) for the input--output relation of the microscale interface solver \(\rmdmix\).
The range of the \(\phaserl\)-phase densities \(\rho_{\cmpa,\phaserl}\), \(\rho_{\cmpb,\phaserl}\) for each component is defined by the convex set, that is formed by the points given in \autoref{tab:mix_data_set_density_corners}.
The barycentric velocity \(v_{\phaser}\) ranges from \(v_{\mathrm{min}} = \SI{-750}{\meter\per\second}\) to
\(v_{\mathrm{max}} = \SI{750}{\meter\per\second}\).
The relative velocities \(v_{\mathrm{rel},\phaserl}\) are bounded by \(v_{\mathrm{rel},\mathrm{min}} = \SI{-500}{\meter\per\second}\) and
\(v_{\mathrm{rel},\mathrm{max}} = \SI{500}{\meter\per\second}\).
Taken all together, the resulting convex set forms the input bounding domain \(B_{\mathrm{in}}\).

\begin{table}[htbp]
    \centering
    \begin{tabular}{
            S[round-mode=places, round-precision=1, table-format=4.1]
            S[round-mode=places, round-precision=1, table-format=4.1]
            S[round-mode=places, round-precision=1, table-format=4.1]
            S[round-mode=places, round-precision=1, table-format=4.1]
        }
        \toprule
        {\(\rho_{\cmpa,\phaser}\) [\si{\kilogram\per\cubic\metre}]} & {\(\rho_{\cmpb,\phaser}\) [\si{\kilogram\per\cubic\metre}]} &
        {\(\rho_{\cmpa,\phasel}\) [\si{\kilogram\per\cubic\metre}]} & {\(\rho_{\cmpb,\phasel}\) [\si{\kilogram\per\cubic\metre}]}                            \\
        \midrule
        1.000000                                                    & 0.000000                                                    & 1024.214739 & 0.000000   \\
        1.000000                                                    & 0.000000                                                    & 0.000000    & 343.403446 \\
        1.000000                                                    & 0.000000                                                    & 1616.252845 & 0.000000   \\
        1.000000                                                    & 0.000000                                                    & 0.000000    & 509.380478 \\
        0.000000                                                    & 1.000000                                                    & 1024.214739 & 0.000000   \\
        0.000000                                                    & 1.000000                                                    & 0.000000    & 343.403446 \\
        0.000000                                                    & 1.000000                                                    & 1616.252845 & 0.000000   \\
        0.000000                                                    & 1.000000                                                    & 0.000000    & 509.380478 \\
        162.996622                                                  & 0.000000                                                    & 1024.214739 & 0.000000   \\
        162.996622                                                  & 0.000000                                                    & 0.000000    & 343.403446 \\
        162.996622                                                  & 0.000000                                                    & 1616.252845 & 0.000000   \\
        162.996622                                                  & 0.000000                                                    & 0.000000    & 509.380478 \\
        0.000000                                                    & 30.960615                                                   & 1024.214739 & 0.000000   \\
        0.000000                                                    & 30.960615                                                   & 0.000000    & 343.403446 \\
        0.000000                                                    & 30.960615                                                   & 1616.252845 & 0.000000   \\
        0.000000                                                    & 30.960615                                                   & 0.000000    & 509.380478 \\
        \bottomrule
    \end{tabular}
    \captionsetup{width=0.75\textwidth}
    \caption{Density corner points for the two-component (argon--methane) model input data set.}
    \label{tab:mix_data_set_density_corners}
\end{table}

\begin{figure}[htbp]
    \centering
    \includegraphics[width=\columnwidth]{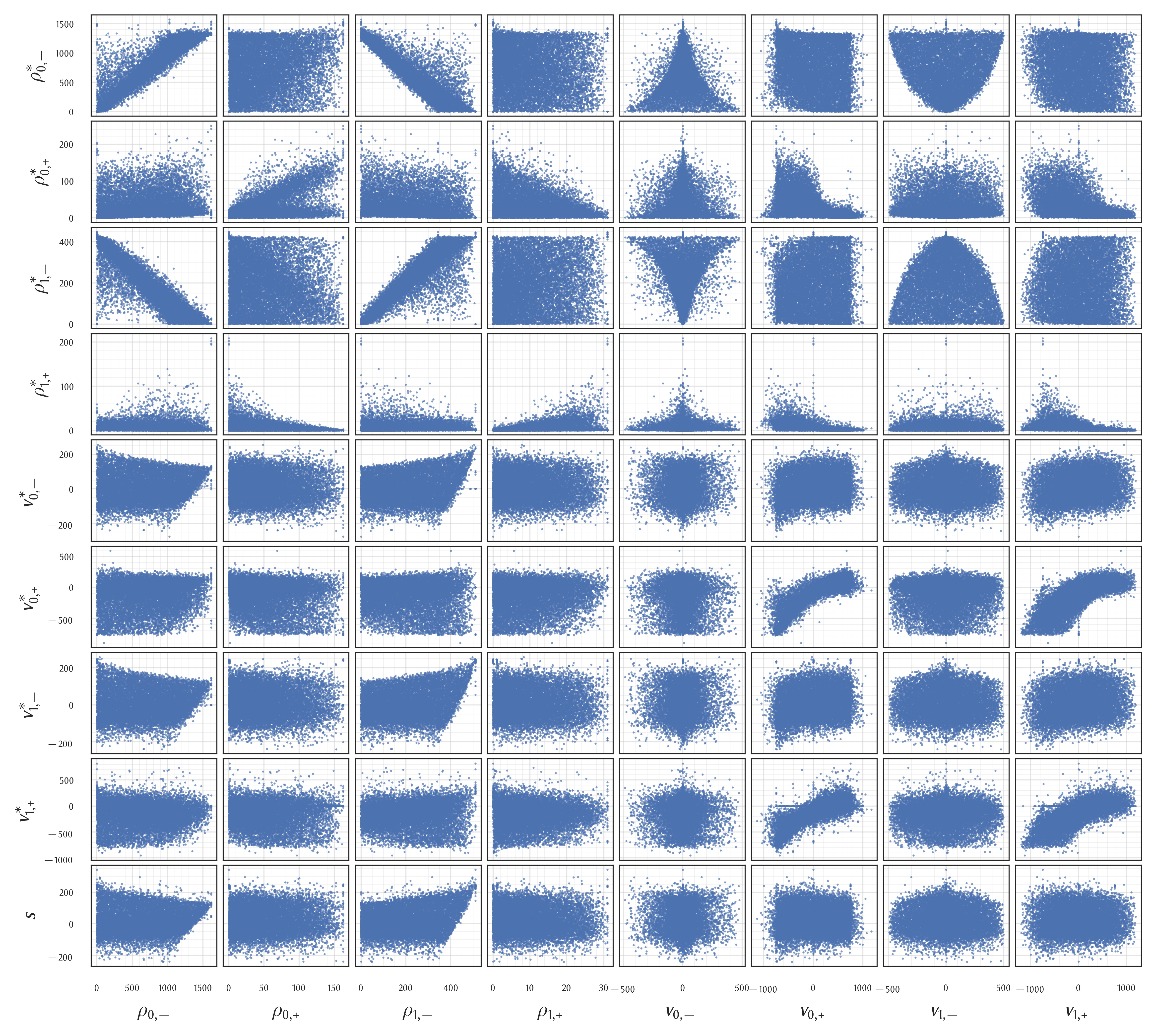}
    \caption{The data set for the isothermal, two-component, two-phase flow multiscale model.
        The data set consists of \num{12000} data points.
        The scatter plots show the relation of the input states
        \(\rho_{\cmpa,\phasel}\), \(\rho_{\cmpb,\phasel}\), \(v_{\cmpa,\phasel}\), \(v_{\cmpb,\phasel}\),
        \(\rho_{\cmpa,\phaser}\), \(\rho_{\cmpb,\phaser}\), \(v_{\cmpa,\phaser}\), \(v_{\cmpb,\phaser}\)
        to their respective wave states
        \(\rho^{*}_{\cmpa,\phasel}\), \(\rho^{*}_{\cmpb,\phasel}\), \(v^{*}_{\cmpa,\phasel}\), \(v^{*}_{\cmpb,\phasel}\),
        \(\rho^{*}_{\cmpa,\phaser}\), \(\rho^{*}_{\cmpb,\phaser}\), \(v^{*}_{\cmpa,\phaser}\), \(v^{*}_{\cmpb,\phaser}\),
        and wave speed \(s\).
        The error-bars show the data range of each separate triplet of MD simulation.}
    \label{fig:mix_dataset:errorbars}
\end{figure}

Using distance-maximizing sampling similar to \cite{mitchell:spectrally:1991},
we generate \(N_{\mathrm{data}} = \num{12000}\) samples in \(B_{\mathrm{in}}\), while exploiting the reduced dimensionality, due to rotational invariance and \(v_{\phasel} \equiv 0\).
Consequently, we generate the input data set
\begin{align*}
    \datasetvar_{\mathrm{in}} = \{ (\rho_{\cmpa,\phasel}, \rho_{\cmpb,\phasel}, v_{\mathrm{rel},\phasel},
    \rho_{\cmpa,\phaser}, \rho_{\cmpb,\phaser}, v_{\phaser}, v_{\mathrm{rel},\phaser})_i \setsep i = 1,\ldots,N_{\mathrm{data}}\}.
\end{align*}
For each \((\rho_{\cmpa,\phasel}, \rho_{\cmpb,\phasel}, v_{\mathrm{rel},\phasel},
\rho_{\cmpa,\phaser}, \rho_{\cmpb,\phaser}, v_{\phaser}, v_{\mathrm{rel},\phaser}) \in \datasetvar_{\mathrm{in}}\)
an MD simulation with the MD interface solver (Algorithm \ref{alg:riemann_mix_euler}) is performed.
It provides with the corresponding input variables the output
\begin{align}
    \label{eq:md:data_evaluation:mix}
    (\rho^{*}_{\cmpa,\phasel}, \rho^{*}_{\cmpb,\phasel},
    m^{*}_{\cmpa,\phasel}, m^{*}_{\cmpb,\phasel},
    \rho^{*}_{\cmpa,\phaser}, \rho^{*}_{\cmpb,\phaser},
    m^{*}_{\cmpa,\phaser}, m^{*}_{\cmpb,\phaser}, s),
\end{align}
which is gathered into the output data set
\begin{align}
    \datasetvar_{\mathrm{out}} = \{ (\rho^{*}_{\cmpa,\phasel}, \rho^{*}_{\cmpb,\phasel},
    m^{*}_{\cmpa,\phasel}, m^{*}_{\cmpb,\phasel},
    \rho^{*}_{\cmpa,\phaser}, \rho^{*}_{\cmpb,\phaser},
    m^{*}_{\cmpa,\phaser}, m^{*}_{\cmpb,\phaser}, s)_i \setsep i = 1,\ldots,N_{\mathrm{data}}\}.
\end{align}
By associating each output data point with its input, we obtain the complete data set
\begin{align*}
    \datasetvar = \Bigl\{
    \bigl(
     & (\rho_{\cmpa,\phasel}, \rho_{\cmpb,\phasel}, v_{\mathrm{rel},\phasel},
    \rho_{\cmpa,\phaser}, \rho_{\cmpb,\phaser}, v_{\phaser}, v_{\mathrm{rel},\phaser}),
    \\
     & (\rho^{*}_{\cmpa,\phasel}, \rho^{*}_{\cmpb,\phasel},
    m^{*}_{\cmpa,\phasel}, m^{*}_{\cmpb,\phasel},
    \rho^{*}_{\cmpa,\phaser}, \rho^{*}_{\cmpb,\phaser},
    m^{*}_{\cmpa,\phaser}, m^{*}_{\cmpb,\phaser}, s)
    \bigr)_i \setsep i = 1,\ldots,N_{\mathrm{data}}
    \Bigr\}.
\end{align*}
Note that we run three MD simulations for each input data point.
The data set \(\datasetvar\) is visualized in %
Figure~\ref{fig:mix_dataset:errorbars}.
It is digitally archived at \cite{magiera:data:2021}.
With the run time of \num{14} to \SI{17}{\minute} for a single, two-component MD simulation, the generation of the whole data set $\datasetvar$ took approximately \SI{3200}{\hour} of computing time, which can be easily split among several of machines to decrease the real time until all data points are sampled.

\subsubsection{Neural Network Training}
\label{sec:network_training:mix}
To train the neural-network surrogate solver \(\networkF_{\nnparams}\), we use a standard training procedure, albeit with some model-specific parameters as described in the following.
The network is comprised of \num{5} hidden layers with \num{60} nodes each.
The data set \(\datasetvar\) from Section \ref{sec:data_generation:mix} is split into a training data set \(\datasetvar_{\mathrm{train}}\) with \num{10800} samples and a validation data set \(\datasetvar_{\mathrm{val}}\) with \num{1200} samples.

Finally, we note that we use a \cres-neural network, which has been developed in \cite{magiera.ray.ea:constraint:2020} to incorporate conservation of selected state variables into the surrogate interface Riemann solver.
For two-component, two-phase flow we decided to preserve the mass conservation across the interface.

\subsection{The Complete HMM\label{sec:HMMfinal}}
Using the machine-learned interface solver \(\networkF_{\nnparams}\) we present the final heterogeneous multiscale method and the interplay of each of its components.
To see the integration of the interface solvers in detail we describe the HMM in the framework of the numerical discretization, a moving-mesh finite-volume method.
The method has been introduced for two space-dimensions in \cite{chalons.rohde.ea:finite:2017} and open-source code for two and three space-dimensions can be found in \cite{alkamper2021interface}.\\
We partition the continuum-scale time interval $[0, t_{\mathrm{end}}]$ according to $0=t^0 < \ldots < t^N = t_{\mathrm{end}}$, and start with an initial mesh ${\mathcal T}(0) = \{ T^0_i \setsep i \in I^0 \}$ for $\contdomain$ that is assumed to consist of simplices $T_i$ indexed via some index set $I^0$.
Let the mesh parameter $h >0$ be given by the maximum of all edges' length.
The mesh ${\mathcal T}(0)$ is supposed to include a connected set of mesh facets of the simplices, which form the approximation of the initial interface $\boldsymbol{\Gamma}_h(0)$, and thus discrete liquid and vapor bulk domains $\boldsymbol{\Omega}_{h,\pm}(0)$.
The initial datum $\UU_0$ from \eqref{eq:iso_euler_mix_iv} is projected to the mesh ${\mathcal T}(0) $ giving
\begin{align} \label{inidata}
    \UU^0_i = \int_{T^0_i} \UU_0(\vx) \, \mathrm{d}\vx \qquad (i\in I^0).
\end{align}
We assume that we have $\UU^0_i \in {\cP}_\pm$ for $T^0_i \in
    \boldsymbol{\Omega}_{h,\pm}(0)$. \\
To introduce the finite-volume time-stepping let us assume that for some discrete time $t^n$, $n \in \{1,\ldots N-1\}$,
we have the same situation, i.e., a simplicial mesh ${\cT}(t^n) = \{ T^n_i \setsep i \in I^n \}$ for $\contdomain$,
a discrete partition $\contdomain = \boldsymbol{\Omega}_{h,-}(t^n) \cup \boldsymbol{\Gamma}_h(t^n) \cup \boldsymbol{\Omega}_{h,+}(t^n)$ with discrete interface
\begin{align}\label{discretegamma}
    \boldsymbol{\Gamma}_h(t^n) = \{ S^n_{ij}\setsep i,j \in I^n, T^n_i \in \boldsymbol{\Omega}_{h,-}(t^n), T^n_j\in \boldsymbol{\Omega}_{h,-}(t^n ) \}
\end{align}
of connected facets, and a set of cell averages $ {\{\UU^n_i\}}_{i\in I^n}$ with the property $\UU^n_i \in {\mathcal P}_\pm$ for $T^n_i \in \boldsymbol{\Omega}_{h,\pm}(t^n)$.
Furthermore, for any facet from $ \boldsymbol{\Gamma}_h(t^n)$ with normal $\vn^n_{ij} \in \bS^{d-1}$ (pointing to the vapor domain $\boldsymbol{\Omega}_{h,+}(t^n)$) we can compute
\[
    (\UU^*_i,\UU^*_j, s^n_{ij} ) = \networkF_{\nnparams}(\UU^n_i,\UU^n_j;\vn^n_{ij} ).
\]
Using the computed speeds $s^n_{ij}$, the moving-mesh finite-volume method from \cite{chalons.rohde.ea:finite:2017}
delivers a mesh
$\widetilde{\cT}(t^{n+1}) = \{ \widetilde{T}^{n+1}_i \setsep i \in \widetilde{I}^{n+1} \}$,
a discrete partition $\contdomain = \widetilde{\boldsymbol{\Omega}}_{h,-}(t^{n+1}) \cup \widetilde{\boldsymbol{\Gamma}}_h(t^{n+1}) \cup \widetilde{\boldsymbol{\Omega}}_{h,+}(t^{n+1}) $, and the family $ {\{\widetilde{\UU}^{n+1}_i\}}_{i\in I^{n+1}}$ with the properties as in time-step $t^n$.
The new mesh $\widetilde{\mathcal T}(t^{n+1})$ evolves from ${\mathcal T}(t^{n+1})$ by affine shifts.
To avoid small cells that would deteriorate the time-step this mesh deformation is followed by a re-meshing leading finally to the new approximations at time $t^{n+1}$, i.e., the mesh
\begin{align*}
    {\mathcal T}(t^{n+1}) = \{ {T}^{n+1}_i \setsep i \in {I}^{n+1} \}, \,
    \contdomain = \boldsymbol{\Omega}_{h,-}(t^{n+1}) \cup \boldsymbol{\Gamma}_h(t^{n+1}) \cup \boldsymbol{\Omega}_{h,-}(t^{n+1}), %
\end{align*}
and the finite-volume approximations
\begin{align*} {\{{\UU}^{n+1}_i\}}_{i\in I^{n+1}}.
\end{align*}
A particular point in the moving-mesh finite-volume method is the choice of the numerical flux functions.
They can be chosen arbitrarily for facets not in $\boldsymbol{\Gamma}_h(t^n)$ but have to be Godunov fluxes across $\boldsymbol{\Gamma}_h(t^n) $ using the adjacent states $\UU^\ast_{i}$, $\UU^\ast_{j}$.
This avoids the occurrence of states not belonging to ${\mathcal P}_{\phaserl}$.
For more details and the control of the explicit time-stepping we refer to \cite{chalons.rohde.ea:finite:2017}.\\
We conclude the section with the complete HMM for two-component, two-phase flow given in algorithmic form.
\begin{boxalgorithm}[label=alg:multiscale]{HMM for Two-Component, Two-Phase Flow}
    \textbf{Input:}
    Initial mesh ${\mathcal T}(0)$ of $\contdomain$ with partition \( \contdomain =\contdomain_{h,\phasel}(0) \cup \continterface_h(0), \contdomain_{h,\phaser}(0)\),
    the family ${\{\UU^0_i\}}_{i\in I^0}$ from \eqref{inidata} with \(\UU^0_i \in {\mathcal P}_\pm\) for $T^0_i \in \contdomain_{h,\phaserl}(0)$.
    \textbf{Prerequisite:}
    the surrogate interface solver \(\networkF_{\nnparams}\), based on \({\cR}_{\mathrm{MD}}\).
    \tcbsubtitle{Algorithm $n \to n+1$}
    \begin{itemize}
        \item For each $S^n_{ij} \in \boldsymbol{\Gamma}_h(t^n)$, see \eqref{discretegamma}, compute
              \begin{align}
                  \label{eq:multiscale:algorithm:interface_solver}
                  (\UU^*_i,\UU^*_j, s^n_{ij} ) = \networkF_{\nnparams} (\UU^n_i,\UU^n_j;\vn^n_{ij} ).
              \end{align}
        \item Run the moving-mesh finite-volume method from \cite{chalons.rohde.ea:finite:2017,alkamper2021interface} with the mesh \(\cT(t^n)\), the partition \( \contdomain =\contdomain_{h,\phasel}(t^n) \cup \continterface_h(t^n) \cup \contdomain_{h,\phaser}(t^n)\), and the family \({\{\UU^n_i\}}_{i \in I^{n+1}}\).
    \end{itemize}
    \textbf{Result:}
    The mesh \(\cT(t^n)\), the partition \( \contdomain =\contdomain_{h,\phasel}(t^{n+1}) \cup \continterface_h(t^{n+1}) \cup \contdomain_{h,\phaser}(t^{n+1})\), and the family \({\{\UU^{n+1}_i\}}_{i \in I^{n+1})}\) with \(\UU^{n+1}_i \in {\mathcal P}_\pm\) for $T^{n+1}_i \in \contdomain_{h,\phaserl}(t^{n+1})$.
\end{boxalgorithm}
Algorithm~\ref{alg:multiscale} provides an approximation \(\UU_h:\contdomain \times [0,t_{\mathrm{end}}] \to \mathcal U\) and an approximation \(\continterface_h(t)\) for the exact interface \(\continterface(t)$,
$ t\in [0,t_{\mathrm{end}}]\), solving the free boundary value problem for \eqref{eq:iso_euler_multicomponent}.
The approximate function $\UU_h$
and the interface $\continterface_h $ are defined by
\[
    \UU_h(\vx,t) = \UU^n_i,
    \quad
    \continterface_h(t) = \continterface_h(t^n) \quad (\vx \in T^n_i, \, t\in [t^n,t^{n+1}]),
\]
using the quantities computed in Algorithm~\ref{alg:multiscale}.
The machine-learned interface solver $\networkF_{\nnparams} $ represents the atomistic microscale whereas the final approximation is valid on the continuum scale only.

\begin{remark}[Discretization of source and non-conservative terms]
    The moving-mesh finite-volume method from \cite{chalons.rohde.ea:finite:2017,MagieraRohde22} uses explicit Euler time-stepping and deals with first-order systems in divergence form only.
    The discretization of the $0$-th-order Maxwell--Stefan diffusion terms is handled by simple evaluation of the cell averages.
    \\
    More complicated is the discretization of the gradients \(\nabla \mu_a\) in the right-hand-side term of the continuum-scale system \eqref{eq:iso_euler_multicomponent}.
    \newline
    In one space-dimension we simply apply central finite-differences in the bulk phases.
    At the cells adjacent to the interface we use left-\slash{}right-sided finite differences, to avoid computing the gradient over the discrete phase boundary.
    \newline
    In two and three space-dimensions, for the $n$-th time-step, we approximate the gradient \(\nabla \mu_a\) in a simplex \(T_j^n\) by linear reconstruction.
    The stencil for the reconstruction includes those neighbors \(T^n_{k}\) of \(T^n_j\) that share a surface with \(T^n_j\) and belong to the same phase domain as \(T^n_j\).
    In that way, we avoid mixing the phases during the reconstruction.
    The reconstruction is done by solving the linear least squares system
    \begin{align}
        \label{eq:chem_pot:lls}
        \mat{A}^\top \mat{A} \vec{\mu}'_a = \mat{A}^\top \vec{b}_a
    \end{align}
    with
    \begin{align*}
        \mat{A}   & =
        \begin{pmatrix}
            \vec{c}_j - \vec{c}_k \\
            \vdots
        \end{pmatrix} \in \bR^{K \times d},
                  &
        \vec{b}_a & =
        \begin{pmatrix}
            \mu_{a,j} - \mu_{a,k} \\
            \vdots
        \end{pmatrix} \in \bR^{K},
                  &
                  & T^n_{k} \text{ in stencil of } T^n_{j},
    \end{align*}
    where \(\vec{c}_j \in \bR^d\) denotes the cell center of \(T^n_j\), \(\mu_{a,j} \in \bR\) the cell value of the chemical potential, and \(K \in \bN\) the number of neighbor cells in the stencil.
    The solution \(\vec{\mu}'_a\) of \eqref{eq:chem_pot:lls} is used as an approximation of the gradient \(\nabla \mu_a\) in cell \(T^n_j\).
\end{remark}
\section{Numerical Simulations for the HMM}
\label{ch:results:mix}

In this section we present a series of numerical results to validate the multiscale method that has been
introduced in Section \ref{ch:multiscale_multicomponent} to solve the free boundary value problem for the two-component, two-phase flow model (\ref{eq:iso_euler_multicomponent}).

The model and numerical parameters used for the simulations in this section, if not otherwise stated, are found in \autoref{sec:parameter_table:mix}.
We choose the Maxwell--Stefan diffusion coefficient to be \(\mathDJ_{\cmpa\cmpb} = 1.0\). We refer to \cite{janzen:diffusion:2019} for the computation of specific Maxwell--Stefan diffusion coefficients leading
to much smaller values.
Numerical tests for simplified settings have shown that the influence of \(\mathDJ_{\cmpa\cmpb}\)
is minor for the overall dynamics, which justifies our choice.
\subsection{One-dimensional Simulation Results}
\label{ch:results:mix:1d}

We consider one-dimensional Riemann problems and compare the multiscale simulation results with MD
simulations on the entire domain as reference.
The section serves as validation for the multi-dimensional droplet simulations in the subsequent Section~\ref{ch:results:mix:2d} and Section~\ref{ch:results:mix:3d}.

\begin{example}[Pressure-driven shock wave]\label{example_1ds}
    In the first example, we simulate a pressure-driven shock wave originating from an argon-rich liquid.
    This is set up by the initial data
    \(\rho_{\cmpa,\phasel} = \SI{1200}{\kilogram\per\cubic\metre}\), \(\rho_{\cmpb,\phasel} = \SI{100}{\kilogram\per\cubic\metre}\) on the liquid side,
    \(\rho_{\cmpa,\phaser} = \SI{10}{\kilogram\per\cubic\metre}\), \(\rho_{\cmpb,\phaser} = \SI{20}{\kilogram\per\cubic\metre}\)
    on the vapor side, and \(\vv_{\cmpa,\phaserl} = \SI{0}{\metre\per\second}\), \(\vv_{\cmpb,\phaserl} = \SI{0}{\metre\per\second}\).
    The simulation results are plotted in \autoref{fig:ar_me_shock}, including the corresponding two-component
    MD simulation results.
    \newline
    In the multiscale solution it can be observed that the liquid phase expands to the right side, while the total density inside the liquid phase decreases.
    Near the interface, argon moves from the liquid into the vapor phase, whereas methane accumulates inside the liquid, and decreases in the vapor phase.
    \newline
    By comparing the multiscale and the MD solution, we see that the interface position is captured accurately.
    Furthermore, on the liquid side, the MD simulation and the multiscale model behave qualitatively in the same manner.
    On the vapor phase side, we see that the argon-component of the multiscale simulation is close to the MD simulation.
    In contrast to that, substantial deviations in the methane-component can be observed.
    We assume that these deviations can be attributed to the isothermal continuum model which does not capture the whole range of the MD simulation.
    The latter is not perfectly isothermal allowing for small temperature variations.
    We observed the same effect for isothermal flow of one component in \cite{MagieraRohde22}.
    Finally, the deviations in the low-density vapor phase might be caused by poor sampling on the MD scale, as few particles of each component are present near the interface.
    This problem can be solved by increasing the number of particles, as well as the size of the computational domain (see Section~\ref{ch:md_mix}) for the MD simulations.
    \begin{figure}[tp]
        \centering
        \includegraphics[width=0.9\columnwidth]{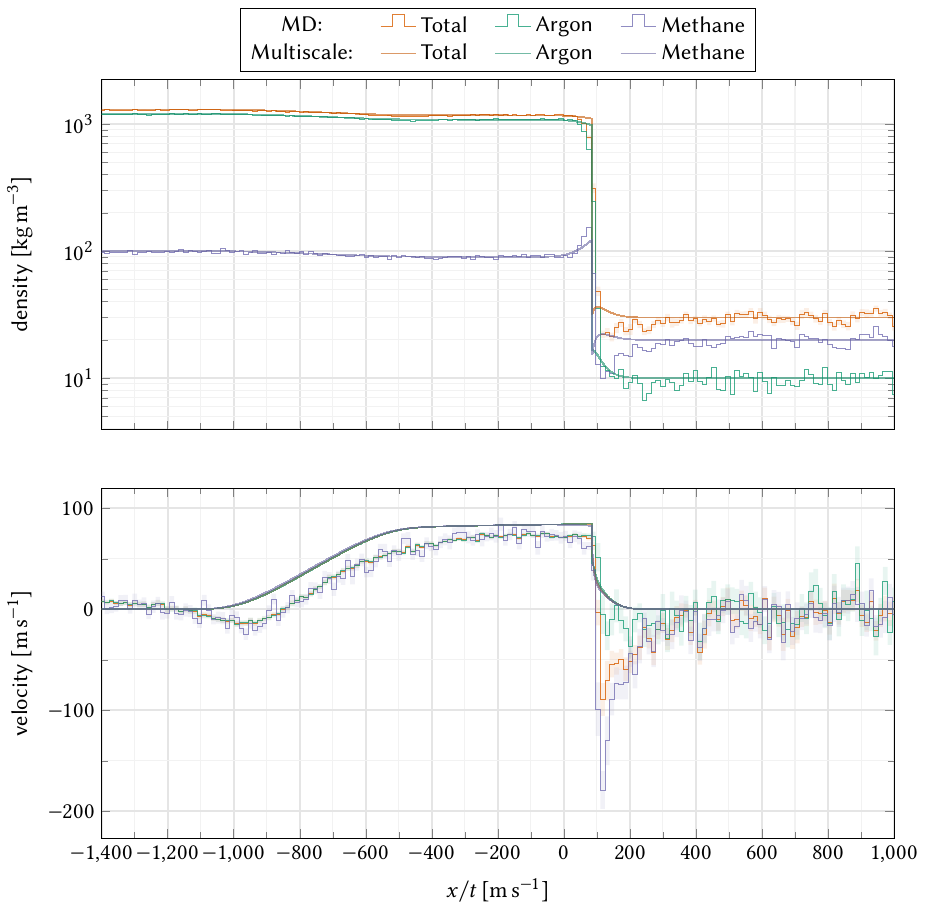}
        \caption{
            Multiscale simulation for the isothermal, two-component, two-phase flow model in one space-dimension, overlaid over the corresponding MD simulations.
        }
        \label{fig:ar_me_shock}
    \end{figure}

\end{example}
\begin{example}[Colliding vapor wave]
    In the second example, we simulate a vapor wave that collides with a liquid argon--methane mixture.
    The corresponding Riemann initial data are
    \(\rho_{\cmpa,\phasel} = \SI{440}{\kilogram\per\cubic\metre}\),
    \(\rho_{\cmpb,\phasel} = \SI{280}{\kilogram\per\cubic\metre}\),
    \(\vv_{\cmpa,\phasel} = \SI{0}{\metre\per\second}\),
    \(\vv_{\cmpb,\phasel} = \SI{0}{\metre\per\second}\)
    for the liquid phase, and
    \(\rho_{\cmpa,\phaser} = \SI{20}{\kilogram\per\cubic\metre}\),
    \(\rho_{\cmpb,\phaser} = \SI{2}{\kilogram\per\cubic\metre}\),
    \(\vv_{\cmpa,\phaser} = \SI{-50}{\metre\per\second}\),
    \(\vv_{\cmpb,\phaser} = \SI{-50}{\metre\per\second}\)
    for the vapor phase.
    The multiscale simulation results are plotted in Figure~\ref{fig:ar_me_condensation}, alongside with their respective MD simulation.
    It can be seen that the vapor wave transmits into the liquid phase, and increases the liquid density slightly.
    Furthermore, the wave speeds, as well as the interface speed, are captured very well by the multiscale model.
    As in Example~\ref{example_1ds} we observe deviations for the velocity in the vapor phase domain.
    \begin{figure}[tp]
        \centering
        \includegraphics[width=0.9\columnwidth]{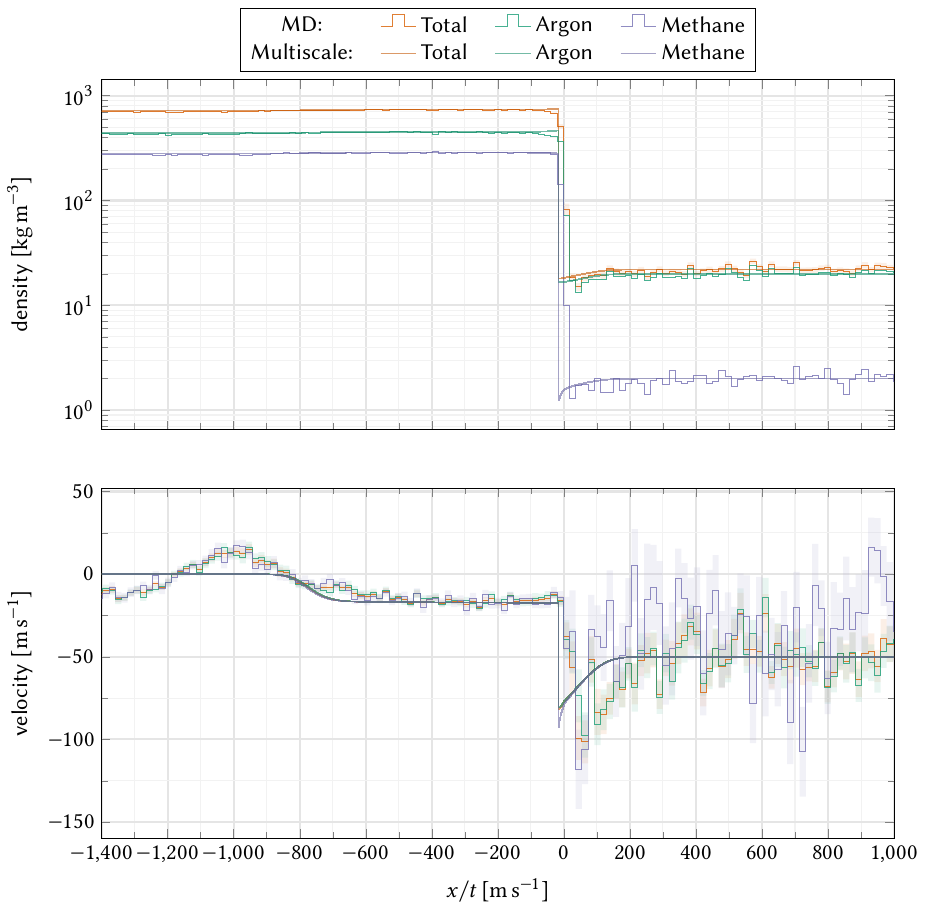}
        \caption{
            Multiscale simulation for the isothermal, two-component, two-phase flow model in one space-dimension, overlaid over the corresponding MD simulations.
            We observe a vapor wave colliding with a liquid argon--methane mixture.
        }
        \label{fig:ar_me_condensation}
    \end{figure}
\end{example}
\subsection{Two-dimensional Simulation Results}
\label{ch:results:mix:2d}

\begin{example}[Interaction shock wave/droplet]\label{example_2ds}
    In this section, we present a two-dimensional, two-component multiscale simulation, where a methane droplet in an argon vapor atmosphere is hit by a shock wave.
    To this end, we consider the domain \(\contdomain = [-1.5, 1.5]^2\), which is split into a liquid droplet
    \(\contdomain_\phasel(0) = \{ \vx \in \bR^2 \setsep \norm{\vx}^2_2 < 0.15\}\)
    and the surrounding vapor domain
    \(\contdomain_\phaser(0) = \contdomain \setminus \overline{\contdomain_\phasel(0)}\).
    The initial conditions are given by
    \begin{align*}
        (\rho_{\cmpa}, \vv_{\cmpa}, \rho_{\cmpb}, \vv_{\cmpb})(\vx, t)
        =
        \begin{cases}
            (180, (0, 0), 400, (0, 0)),  & \vx \in \contdomain_\phasel(0),                          \\
            (20, (0, 0), 4, (0, 0)),     & \vx \in \contdomain_\phaser(0) \text{ and } x \geq -0.5, \\
            (20, (150, 0), 4, (150, 0)), & \vx \in \contdomain_\phaser(0) \text{ and } x < -0.5.
        \end{cases}
    \end{align*}
    On the left side, at \(x = -1.5\) we apply inflow boundary conditions by setting the ghost cell value to
    $(\rho_{\cmpa},\allowbreak \vv_{\cmpa},\allowbreak \rho_{\cmpb},\allowbreak \vv_{\cmpb})
        = (20, \allowbreak (150, \allowbreak 0), \allowbreak 4, \allowbreak (150,\allowbreak 0))$.
    On the opposing side, at \(x = 1.5\), we use outflow boundary conditions.
    At the top and bottom, i.e., at \(y = \pm 1.5\), reflecting boundary conditions simulate closed walls.

    The simulation results are depicted in
    \autoref{fig:mix:density:2d} and
    \autoref{fig:mix:mole_fraction_0:2d}.
    In the first figure, the component-wise densities and velocities are shown, and in the second figure the argon-mole fraction is plotted.
    \newline
    In the beginning, the liquid and vapor phase are not in equilibrium, resulting in small oscillations of the droplet.
    The oscillations can be clearly observed in the time evolution of the \(\phasel\)-averaged quantities, as illustrated in Figure~\ref{fig:2d_mix_data_averages}.
    Then the vapor wave hits the droplet, causing a ripple moving through its surface and finally pushing it through the vapor atmosphere.
    Throughout the simulation, argon accumulates inside the liquid phase, leading to a growth of the droplet -- see Figure~\ref{fig:2d_mix_data_averages}.
    \begin{figure}[tp]
        \centering
        \includegraphics[width=0.75\columnwidth]{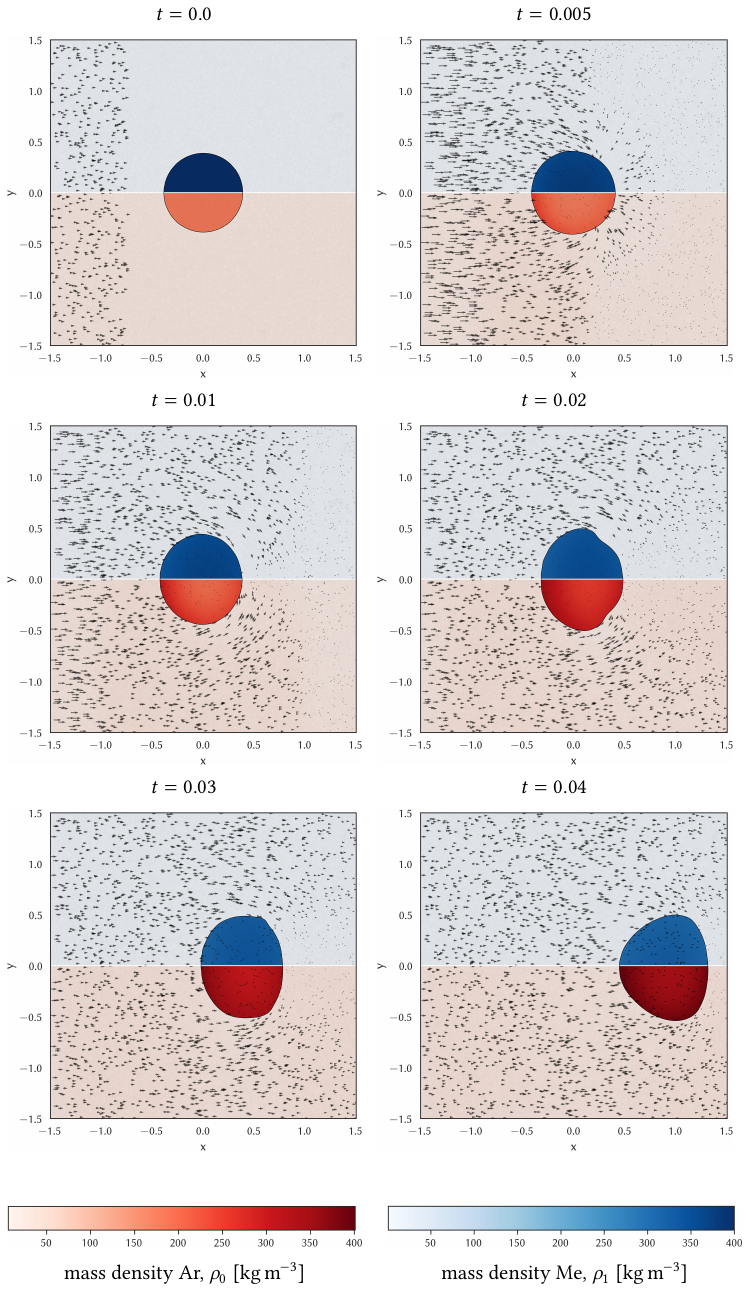}
        \caption{
            Two-dimensional multiscale simulations.
            Each sub-figure depicts the densities \(\rho_{\cmpa}\), \(\rho_{\cmpb}\) and velocities \(\vv_{\cmpa}\), \(\vv_{\cmpb}\) of each component at various time steps.
            The upper part of each sub-figure shows \(\rho_{\cmpb}\), \(\vv_{\cmpb}\) for methane, and the lower part \(\rho_{\cmpa}\), \(\vv_{\cmpa}\) for argon.
        }
        \label{fig:mix:density:2d}
    \end{figure}
    \begin{figure}[tp]
        \centering
        \includegraphics[width=0.8\columnwidth]{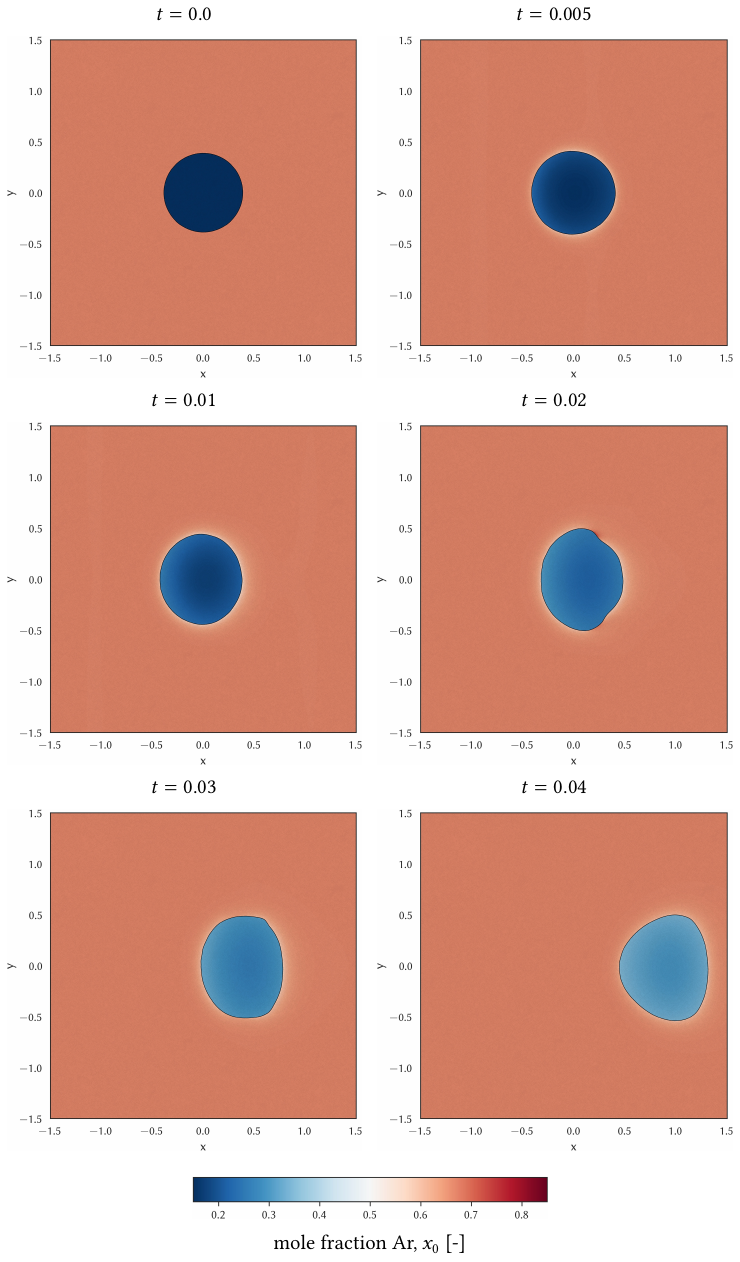}
        \caption{
            Two-dimensional multiscale simulation.
            Each sub-figure depicts the mole fraction \(x_{\cmpa}\) of argon at various time steps.
        }
        \label{fig:mix:mole_fraction_0:2d}
    \end{figure}
    \begin{figure}[htbp]
        \centering
        \includegraphics[width=0.82\columnwidth]{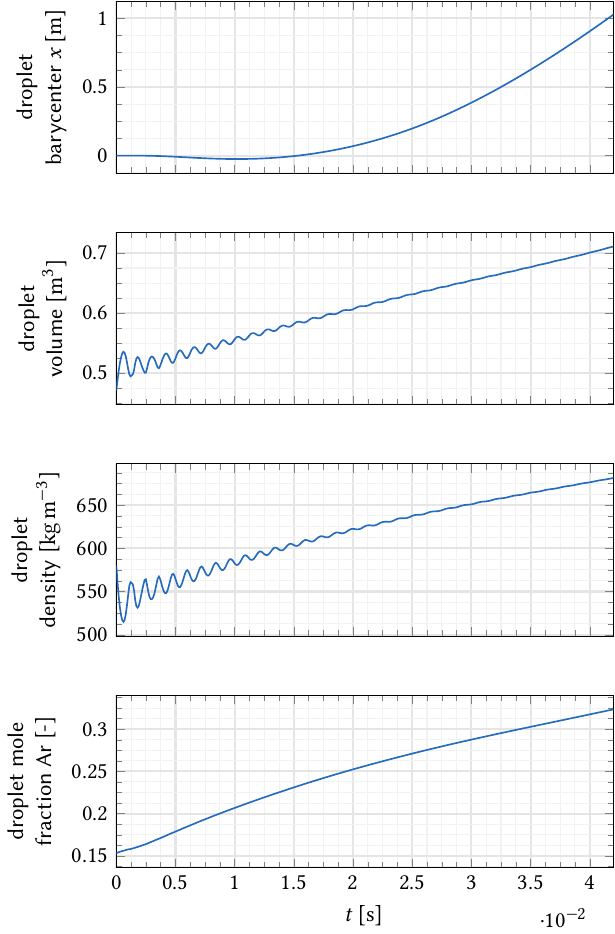}
        \caption{
            Time evolution of liquid phase-averaged quantities for the two-component two-phase flow multiscale simulation in two space-dimensions.
        }
        \label{fig:2d_mix_data_averages}
    \end{figure}

\end{example}
\subsection{Three-dimensional Simulation Results}
\label{ch:results:mix:3d}
\begin{example}[Interaction shock wave/droplet]\label{example_3ds}
    The last example in this section is a three-dimensional version of Example~\ref{example_2ds}, i.e., a simulation of a liquid droplet, that consists mostly of methane, inside an argon--methane vapor atmosphere.
    The domain \(\contdomain = [-5, 5]^3\) is split into a liquid droplet
    \(\contdomain_\phasel(0) = \{ \vx \in \bR^3 \setsep \norm{\vx}_2 < 2\}\)
    and the surrounding vapor domain
    \(\contdomain_\phaser(0) = \contdomain \setminus \overline{\contdomain_\phasel(0)}\).
    For the initial data, we set
    \begin{align*}
        (\rho_{\cmpa}, \vv_{\cmpa}, \rho_{\cmpb}, \vv_{\cmpb})(\vx, t)
        =
        \begin{cases}
            (180, (0, 0, 0), 400, (0, 0, 0)),  & \vx \in \contdomain_\phasel(0),                        \\
            (20, (0, 0, 0), 4, (0, 0, 0)),     & \vx \in \contdomain_\phaser(0) \text{ and } x \geq -3, \\
            (20, (200, 0, 0), 4, (200, 0, 0)), & \vx \in \contdomain_\phaser(0) \text{ and } x < -3.
        \end{cases}
    \end{align*}
    Inflow boundary conditions are applied at \(x = -5\), by setting the ghost cell values to
    \((\rho_{\cmpa},\allowbreak \vv_{\cmpa},\allowbreak \rho_{\cmpb},\allowbreak \vv_{\cmpb})
    = (20,\allowbreak (200,\allowbreak 0,\allowbreak 0), \allowbreak 4, \allowbreak (200,\allowbreak 0,\allowbreak 0))\).
    Opposed to that side, at \(x = 5\), outflow boundary conditions are applied.
    On every other side of the domain we implement reflecting boundary conditions.

    For the multiscale simulation results we refer to Figure~\ref{fig:mix:shock_droplet:density:3d} and Figure~\ref{fig:mix:shock_droplet:interface:3d}.
    Both figures show the three-dimensional solution visualized on the plane through \(z=0\).
    In Figure~\ref{fig:mix:shock_droplet:density:3d} the component-wise densities and velocities are depicted.
    We display also the (projected) mesh to demonstrate that the interface is resolved on the discrete level.
    In Figure~\ref{fig:mix:shock_droplet:interface:3d} the phase boundary is shown alongside the barycentric velocity magnitude.
    \newline
    We observe that the droplet is hit by a shock wave, which results in a ripple that runs over the liquid surface.
    Additionally, the momentum from the vapor causes the droplet to move to the right side.
    Considering the two fluid components, we see that the methane concentration inside the liquid droplet decreases slowly, whereas the amount of liquid argon slowly increases.

    Next, we provide an overview over the simulation run times.
    A single time step of the simulation takes on average
    \SI[round-mode = places, round-precision = 1]{2.70159662323362084786}{\second}.
    This splits approximately into \SI{1.0}{\second} for the moving mesh operations (see Section \ref{sec:HMMfinal}), and \SI{1.7}{\second} for the finite volume part of the simulation.
    The mesh consists of
    \num[scientific-notation=true, round-mode = places, round-precision = 1]{273392.22}
    cells, with
    \num[scientific-notation=true, round-mode=places, round-precision=1]{2043.442080}
    interface facets
    averaged over all \num{100000} time steps.
    In total the whole simulation takes
    \SI[round-mode=places, round-precision=0]{75.044}{\hour} on a single desktop computer.
    \newline
    Note that in this simulation the need for a surrogate solver becomes apparent.
    If we would not have employed a surrogate, we would need to perform an MD simulation (taking \SIrange{14}{17}{\minute}) for every interface facet at every single time step.
    For a single time step, that would accumulate to approximately
    \SI{500}{\hour}
    of computational time (albeit parallelizable).
    Compared to that, the surrogate solver needs \SI{0.102}{\milli\second} for a single evaluation.
    The offline phase of the surrogate solver takes approximately \SI{3200}{\hour} plus training.
    The majority of this workload (generating the training data set) is however embarrassingly parallel.
    Consequently, even if we include the offline phase, using a surrogate solver pays off after a few time steps.
    \begin{figure}[tbp]
        \centering
        \includegraphics[width=0.76\columnwidth]{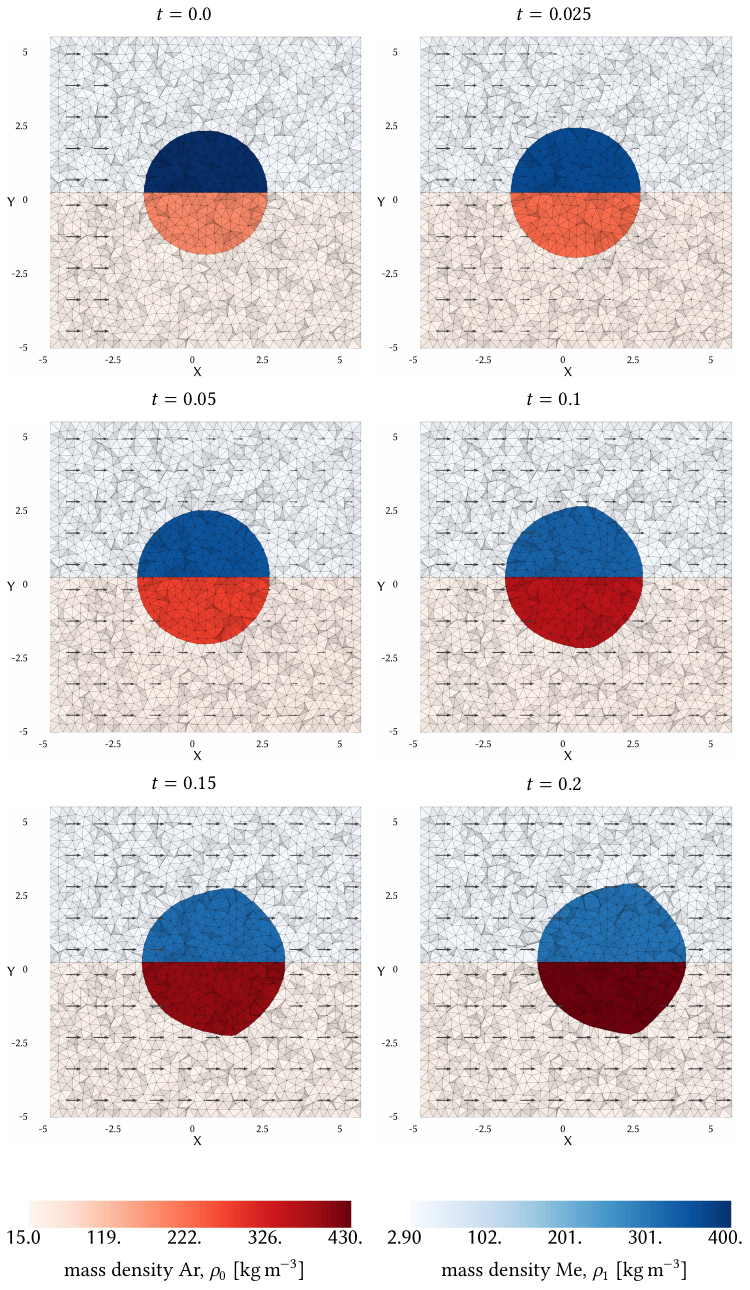}
        \caption{
            Three-dimensional multiscale simulation of the two-component two-phase flow model.
            The upper part of each sub-figure shows \(\rho_{\cmpb}\), \(\vv_{\cmpb}\) for methane, and the lower part \(\rho_{\cmpa}\), \(\vv_{\cmpa}\) for argon, at various time steps.}
        \label{fig:mix:shock_droplet:density:3d}
    \end{figure}
    \begin{figure}[tbp]
        \centering
        \includegraphics[width=0.8\columnwidth]{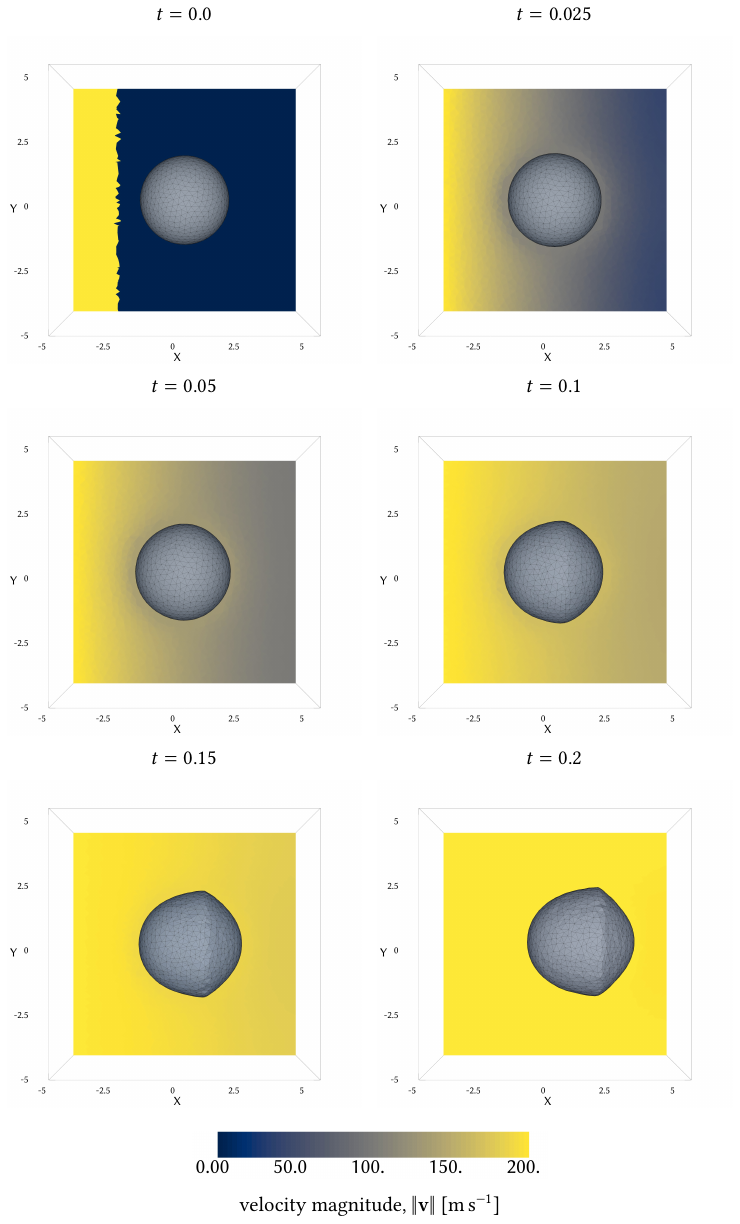}
        \caption{
            Three-dimensional multiscale simulation.
            Each sub-figure depicts the barycentric velocity magnitude \(\norm{\vv}_2\) and the phase boundary at various time steps.}
        \label{fig:mix:shock_droplet:interface:3d}
    \end{figure}

\end{example}

\section{Conclusions}\label{sec:summary}
The main result of our study is a proof-of-concept validation for a HMM that has been designed to
model the dynamics of compressible two-component flow with liquid--vapor phase transition.
The behavior of the numerical solutions on the continuum scale matches the MD simulations, with even quantitatively correct results in the liquid phase domains.
Improving the results in the vapor phase seems to require more computational effort in the sampling procedure and the thermalization process.
It is noteworthy that we do not consider a generic fluid mixture, but use a real EOS for argon--methane mixtures.
Another significant outcome is that the HMM is able to simulate two- and three-dimensional droplets in case of complex interactions.
At the same time, both small and large-scale deformations of the (sharp) phase boundary can be rendered by the moving-mesh finite-volume method during the numerical simulations. \\
Up to our knowledge, this is the first time that compressible mixtures of real fluids with resolved sharp interfaces are simulated using an MD-based and neural-network accelerated multiscale model.
In particular, we are able to simulate phase transitions of fluid mixtures -- without the need to prescribe some (ad-hoc) closure relations at the interface.\\
We think that the HMM in combination with an ML surrogate modelling paves the way to handle more complex flow scenarios.
Surface tension has not been taken into account, despite its importance for two-phase evolution.
Here, it is possible to account for curvature effects on the continuum-scale by adding a geometric term to the momentum
equation \cite{rohde.zeiler:riemann:2018}.
However, the surface tension coefficient might then require additional MD simulations like for the EOS.
On the expense of more MD simulations we expect that the extension to mixtures with more than two components is straightforward. A conceptual advantage of the underlying moving-mesh finite-volume method is that tailored numerical schemes can be applied in the liquid and the vapor phase.
This would allow for benefitting from e.g. implicit time-stepping in the low-Mach liquid domain and -- more important for mixtures -- to account for large Maxwell--Stefan coefficients.
\medskip

{\bfseries{Acknowledgments.}}
This work was supported by the Deutsche Forschungsgemeinschaft (DFG, German Research Foundation) through the project \mbox{SFB--TRR 75} with the project number 84292822,
and the DFG under Germany's Excellence Strategy - EXC 2075 with the project number 390740016.

\printbibliography
\clearpage
\appendix
\section{Parameter Tables}\label{sec:parameter_table:mix}
\begin{table}[ht!]
    \centering \begin{tabular}{lr}
        \toprule
        \multicolumn{2}{c}{{Molecular Dynamics}}       \\ \midrule
        cutoff-radius
        $r_{\mathrm{cutoff}}$            & \num{2.5}   \\

        total number of particles
        $N$                              & \num{32768} \\
        number of time steps
        $N_{\mathrm{end}}$               & \num{5e2}   \\
        time step
        $\Delta \tau$                    & \num{5e-4}  \\
        time sampling ratio
        $\fractionvar_{\mathrm{t-smpl}}$ & \num{0.2}   \\
        \bottomrule                                    \\
    \end{tabular}
    \captionsetup{width=0.75\textwidth}
    \caption{Parameter table for the atomistic interface solver.}
    \label{tab:mix_paramater}
\end{table}
\begin{table}[ht!]
    \centering
    \begin{tabular}{lr}
        \toprule
        \multicolumn{2}{c}{\phantom{wwwwwww}Simulations for $d=1$\phantom{wwwwwww} } \\ \midrule
        end time $t_{\mathrm{end}}$ & \num{0.003}                                    \\
        time step $\Delta t$        & \num{e-7}                                      \\
        domain $\contdomain$        & $[-5, 5]$                                      \\
        base cell width $\Delta x$  & \num{2e-3}                                     \\
        $\Delta x_{\mathrm{min}}$   & $\sfrac{1}{2} \Delta x$                        \\
        $\Delta x_{\mathrm{max}}$   & $\sfrac{3}{2} \Delta x$                        \\
        Maxwell--Stefan diffusion   &                                                \\ coefficient $\mathDJ_{\cmpa\cmpb}$ & \num{1.0} \\
        \bottomrule
    \end{tabular}
    \captionsetup{width=0.75\textwidth}
    \caption{Parameters for Section \ref{ch:results:mix:1d}.}
    \label{tab:mix_paramater:mmfv:1d}
\end{table}
\begin{table}[ht]
    \centering
    \begin{tabular}{lr}
        \toprule
        \multicolumn{2}{c}{\phantom{wwwwwww}Simulations for $d=2$\phantom{wwwwwww}} \\ \midrule
        time step $\Delta t$        & \num{5e-7}                                    \\
        domain $\contdomain$        & $[-1.5, 1.5]^2$                               \\
        base edge length $\Delta x$ & \num{1.25e-2}                                 \\
        $\Delta x_{\mathrm{min}}$   & $\sfrac{1}{2} \Delta x$                       \\
        $\Delta x_{\mathrm{max}}$   & $\sfrac{3}{2} \Delta x$                       \\
        \bottomrule
    \end{tabular}
    \captionsetup{width=0.75\textwidth}
    \caption{Parameters for Section \ref{ch:results:mix:2d}.}
    \label{tab:mix_paramater:mmfv:2d}
\end{table}
\begin{table}[ht]
    \centering
    \begin{tabular}{lr}
        \toprule
        \multicolumn{2}{c}{\phantom{wwwwwww}Simulations for $d=3$\phantom{wwwwwww}} \\ \midrule
        time step $\Delta t$        & \num{2e-6}                                    \\
        domain $\contdomain$        & $[-5, 5]^3$                                   \\
        base edge length $\Delta x$ & \num{0.25}                                    \\
        $\Delta x_{\mathrm{min}}$   & $\sfrac{1}{2} \Delta x$                       \\
        $\Delta x_{\mathrm{max}}$   & $\sfrac{3}{2} \Delta x$                       \\
        Maxwell--Stefan diffusion   &                                               \\ coefficient $\mathDJ_{\cmpa\cmpb}$ & \num{1.0} \\
        \bottomrule
    \end{tabular}
    \captionsetup{width=0.75\textwidth}
    \caption{Parameters for Section \ref{ch:results:mix:3d}.}
    \label{tab:mix_paramater:mmfv:3d}
\end{table}
\end{document}